\documentclass[12pt]{article}
\usepackage{amsmath,amsfonts}
\usepackage{graphicx, xcolor}
\usepackage{enumerate, bm, dsfont}
\usepackage{natbib,microtype,subfig, booktabs}
\usepackage[top=1.3in, bottom=1.3in, left=1.3in, right=1.3in]{geometry}
\usepackage{url} 

\newcommand{\blind}{1}
\usepackage[pdftex,colorlinks=true]{hyperref}
\definecolor{darkblue}{rgb}{0,0,.6}
\hypersetup{citecolor=darkblue,linkcolor=darkblue,urlcolor=darkblue}

\newcommand{\X}{\mathcal{X}}
\graphicspath{{plots/}}
\DeclareMathOperator*{\argmin}{\arg\!\min}

\newtheorem{theorem}{Theorem}[section]
\newtheorem{lem}{Lemma}[section]
\newtheorem{corollary}{Corollary}[section]
\newtheorem{prop}{Proposition}[section]

\date{}

\begin{document}

\def\spacingset#1{\renewcommand{\baselinestretch}%
{#1}\small\normalsize} \spacingset{1}


\if1\blind
{
  \title{\bf Bootstrap prediction bands for functional time series}
  \author{Efstathios Paparoditis
  \hspace{.2cm}\\
    Department of Mathematics and Statistics\\ University of Cyprus\\
    and \\
    Han Lin Shang \\
    Department of Actuarial Studies and Business Analytics\\ Macquarie University}
  \maketitle
} \fi

\if0\blind
{
  \bigskip
  \bigskip
  \bigskip
  \begin{center}
    {\LARGE\bf Bootstrap prediction bands for functional time series}
\end{center}
  \medskip
} \fi

\bigskip
\begin{abstract}
A bootstrap procedure for constructing prediction bands for a stationary functional time series is proposed. The procedure exploits a general vector autoregressive representation of the time-reversed series of Fourier coefficients appearing in the Karhunen-Lo\`{e}ve representation of the functional process. It generates backward-in-time functional replicates that adequately mimic the dependence structure of the underlying process in a model-free way and have the same conditionally fixed curves at the end of each functional pseudo-time series. The bootstrap prediction error distribution is then calculated as the difference between the model-free, bootstrap-generated future functional observations and the functional forecasts obtained from the model used for prediction. This allows the estimated prediction error distribution to account for the innovation and estimation errors associated with prediction and the possible errors due to model misspecification. We establish the asymptotic validity of the bootstrap procedure in estimating the conditional prediction error distribution of interest, and we also show that the procedure enables the construction of prediction bands that achieve (asymptotically) the desired coverage. Prediction bands based on a consistent estimation of the conditional distribution of the studentized prediction error process also are introduced. Such bands allow for taking more appropriately into account the local uncertainty of the prediction. Through a simulation study and the analysis of two data sets, we demonstrate the capabilities and the good finite-sample performance of the proposed method.
\end{abstract}

\noindent%
\vfill

\newpage
\spacingset{1.5} 

\section{Introduction}\label{sec:intro}

Functional time series consist of random functions observed at regular time intervals. They can be classified into two main categories depending on whether the continuum is also a time variable. First, functional time series can arise from measurements obtained by separating an almost continuous time record into consecutive intervals \citep[e.g., days, weeks, or years; see, e.g.,][]{HK12}. We refer to such data structures as sliced functional time series, examples of which include daily price curves of a financial stock \citep{KRS17} and intraday particulate matter \citep{Shang17}. When the continuum is not a time variable, functional time series can also arise when observations over a period are considered as finite-dimensional realizations of an underlying continuous function \citep[e.g., yearly age-specific mortality rates; see, e.g.,][]{CM09, HS09}.

In either case, the underlying strictly stationary stochastic process is denoted by ${\bf X}= \{\X_t, t \in \mathbb{Z}\}$,  where each $\X_t$ is a random element in a separable Hilbert space ${\mathcal H}$, with values $\X_t(\tau)$ and $\tau$ varying within a compact interval $ \mathcal{I}\subset \mathbb{R}$. We assume that ${\mathcal I}=[0,1]$ without loss of generality. Furthermore, we assume that $ {\mathcal H}=L^2$,  the set of (equivalence classes of) measurable, real-valued functions $x(\cdot)$ defined on $[0,1]$ and satisfying $\int_0^1x^2(\tau)d\tau <\infty$. Central statistical issues include modeling of the temporal dependence of the functional random variables $ \{\X_t, t\in \mathbb{Z}\}$, making inferences about parameters of interest, and predicting future values of the process when an observed stretch $\X_1, \X_2, \ldots, \X_n$ is given. Not only is it vital to obtain consistent estimators, but to also estimate the uncertainty associated with such estimators, the construction of confidence or prediction intervals, and the implementation of hypothesis tests \citep[e.g.,][]{HKR14}. When such inference problems arise in functional time series, a resampling methodology, especially bootstrapping, is an important alternative to standard asymptotic considerations. For independent and identically distributed (i.i.d.) functional data, bootstrap procedures have been considered, among others, by \cite{CFF06, MP11, GGC13, Shang15, PS16}, where appropriate sampling from the observed sample is used to mimic sampling from the population. However, for functional time series, the existing temporal dependence between the random elements $\X_t$ significantly complicates matters, and the bootstrap must be appropriately adapted in order to be successful.

The development of bootstrap procedures for functional time series has received increasing attention in recent decades. In an early paper, \cite{PR94} obtained weak convergence results for approximate sums of weakly dependent, Hilbert space-valued random variables. \cite{DSW15} also obtained weak convergence results for Hilbert space valued random variables, which are assumed to be weakly dependent in the sense of near-epoch dependence and showed consistency of a non-overlapping block bootstrap procedure. \cite{RAV15} extended the stationary bootstrap to functional time series, \cite{FV11} applied a residual-based bootstrap procedure to construct confidence intervals for the regression function in a nonparametric setting, and \cite{ZhuPoli17} to kernel estimators. \cite{FN16} proposed a residual-based bootstrap procedure for functional autoregressions. \cite{PPS18} established theoretical results for the moving block and the tapered block bootstrap, \cite{Shang18} applied a maximum entropy bootstrap procedure, and \cite{Paparoditis18} proposed a sieve bootstrap  for functional time series.   

In this paper, we build on the developments mentioned above and focus on constructing prediction intervals or bands for a functional time series. To elaborate, suppose that for every $t\in \mathbb{Z}$, the zero mean random element $\X_t$ is generated as 
\begin{equation} \label{eq.process}
\X_t=f(\X_{t-1}, \X_{t-2}, \ldots) + \varepsilon_t,
\end{equation}
where  $ f : {\mathcal H}^\infty \rightarrow {\mathcal H} $  is some appropriate operator and  $ \{\varepsilon_t\}$ is  a zero mean  i.i.d. innovation process in ${\mathcal H}$ with  $ E\|\varepsilon_t\|^2 <\infty$ and $ \|\cdot\|$ the norm of ${\mathcal H}$. For simplicity, we write $ \varepsilon_t \sim i.i.d.(0,C_\varepsilon)$, where $C_\varepsilon=E (\varepsilon_t \otimes\varepsilon_t )$ is the covariance operator of $ \varepsilon_t$ and $ \otimes$ denotes  the tensor operator,  defined by  $ (x  \otimes y )(\cdot)= \langle x,\cdot  \rangle y $,  for $ x,y \in {\mathcal H}$. Suppose that based on the last $k$ observed functional elements,  $\X_n, \X_{n-1}, \ldots, \X_{n-k+1}$,  $ k<n$, a predictor 
\begin{equation} \label{eq.GenPred}
  \widehat{\X}_{n+h}= \widehat{g}_{(h)}( \X_n, \X_{n-1}, \ldots, \X_{n-k+1}),
\end{equation}
of $ \X_{n+h}$ is used, where $h\in \mathbb{N}$ is the prediction horizon and $ \widehat{g}_{(h)} :{\mathcal H}^k\rightarrow {\mathcal H}$ some estimated operator. For instance, we may think of $ \widehat{g}_{(h)}$ as an estimator of the best predictor, i.e., of the conditional expectation $ E(\X_{n+h}|\X_n,\X_{n-1}, \ldots, \X_{n-k+1})$. The important case, however, we have in mind,  is that where a model 
\begin{equation} 
\X_t=g(\X_{t-1}, \X_{t-2}, \ldots, \X_{t-k}) + v_t,  \label{eq:model}
\end{equation}
is used to predict $ \X_{n+h}$. Here  $g: {\mathcal H}^k \rightarrow {\mathcal H} $  is  an unknown, linear bounded operator and $v_t \sim i.i.d. (0, C_v)$.  Using model~\eqref{eq:model},  a $h$-step-ahead predictor can be obtained as
\begin{equation} \label{eq.pred1}
\widehat{\X}_{n+h} = \widehat{g}\left(\widehat{\X}_{n+h-1}, \widehat{\X}_{n+h-2}, \ldots, \widehat{\X}_{n+h-k}\right),
\end{equation}
where $ \widehat{g}$ is an estimator of the operator $g$ in~\eqref{eq:model} and $ \widehat{\X}_t \equiv \X_t$, if $ t \in \{n-k+1, n-k+2, \ldots, n\}$. Notice that  the predictor~\eqref{eq.pred1} also can be written as $\widehat{\X}_{n+h} = \widehat{g}_{(h)}(\X_n,  \X_{n-1},\ldots, \X_{n-k+1})$ for some appropriate operator
$ \widehat{g}_{(h)} $. In particular, setting $\widehat{g}_{(1)}=\widehat{g}$,  it is easily seen that  
$ \widehat{g}_{(h)}=\widehat{g}(\widehat{g}_{(h-1)}, \ldots, \widehat{g}_{(1)}, \X_n, \ldots, \X_{n-k+h})$ for $ 2   \leq h \leq k$ and 
$ \widehat{g}_{(h)} =\widehat{g}(\widehat{g}_{(h-1)}, \ldots, \widehat{g}_{(h-k)})$ for $ h>k$.
We stress here the fact that we do not assume that model~\eqref{eq:model} used for prediction coincides with the true data generating process (\ref{eq.process}); that is, we allow for model misspecification. Observe that the simple case, where $g$  is known up to a finite-dimensional vector of parameters, is also covered by the above setup.

As already mentioned, our aim is to  construct a prediction  band for $ \X_{n+h}$ associated with the predictor $ \widehat{\X}_{n+h}$. That is, given the  part $ {\mathcal X}_{n,k} = (\X_n, \X_{n-1}, \ldots, \X_{n-k+1})$ of the functional time series observed and for any $\alpha \in (0,1)$,  we want to construct a band  denoted by  $\{ [\widehat{\X}_{n+h}(\tau)-L_{n,h}(\tau), \ \widehat{\X}_{n+h}(\tau)+ U_{n,h}(\tau)], \ \tau \in [0,1]\}$ such that 
\[ \lim_{n\rightarrow \infty} P\Big( \widehat{\X}_{n+h}(\tau)-L_{n,h}(\tau) \leq  \X_{n+h}(\tau) \leq  \widehat{\X}_{n+h}(\tau)+ U_{n,h}(\tau), \ \mbox{for all} \ \tau \in [0,1] \Big| \X_{n,k} \Big)=1-\alpha.\] Toward this end we focus on the estimation of the conditional distribution of  the prediction  error  $ \mathcal E_{n+h}=  \X_{n+h} - \widehat{\X}_{n+h}$ given $ {\mathcal X}_{n,k}$, which is a key quantity  for the construction of the prediction band of interest. Using~\eqref{eq.pred1}, this error   can  be decomposed as 
\begin{align*}
\mathcal{E}_{n+h} & :=  \X_{n+h} - \widehat{\X}_{n+h}\\
 & = \varepsilon_{n+h}\\
 & \ \  \ + \Big[f\left(\X_{n+h-1}, \X_{n+h-2}, \ldots\right) - g\left(\X_{n+h-1}, \X_{n+h-2}, \ldots, \X_{n+h-k}\right) \Big]\\
 & \ \  \ + \Big[g\left(\X_{n+h-1}, \X_{n+h-2}, \ldots, \X_{n+h-k}\right) - \widehat{g}\left(\widehat{\X}_{n+h-1}, \widehat{\X}_{n+h-2}, \ldots, \widehat{\X}_{n+h-k}\right)\Big]\\
 & = \mathcal{E}_{I,n+h}+\mathcal{E}_{M,n+h} + \mathcal{E}_{E,n+h},
\end{align*}
with an obvious notation for $\mathcal{E}_{I,n+h}$, $\mathcal{E}_{M,n+h}$ and $\mathcal{E}_{E,n+h}$. Notice that  $\mathcal{E}_{I,n+h}$ is the error attributable to the i.i.d. innovation, $\mathcal{E}_{M,n+h} $ is the model specification error, and $\mathcal{E}_{E,n+h}$ is the error attributable to estimation of the unknown operator $g$ and of the random elements $ \X_{n+h-1}, \ldots, \X_{n+h-k}$ used for $h$-step prediction.  Observe that if  $ h=1$ then $\mathcal{E}_{E,n+1}$ only depends on the estimation error $ \widehat{g}-g$.  Furthermore, if $ \widehat{g}$ is a consistent estimator of $g$, for instance, if $ \|\widehat{g}-g\|_{\mathcal L} \stackrel{P}{\rightarrow} 0$, with $ \|\cdot\|_{\mathcal L}$ being the operator norm, the part of the estimation error $\mathcal{E}_{E,n+h} $ which is due to the estimator $ \widehat{g}$  is  asymptotically negligible. On the contrary, the misspecification error $ \mathcal{E}_{M,n+h}$ may not vanish asymptotically if the model used for prediction is different from the one generating the data.  

To better illustrate the above discussion, consider the following example. Suppose that $ \X_t$ is generated according to the  FAR(2) model $\X_t=\Phi_1(\X_{t-1}) + \Phi_2(\X_{t-2}) +\varepsilon_t$, $ \Phi_2 \neq 0$,  and that a FAR(1) model $ \X_t=R(\X_{t-1}) +v_t$ is used for prediction, where  $ \Phi_1$, $\Phi_2$, and $ R$ are appropriate operators.  In general, $ R \neq \Phi_1$. With  $ \widehat{R}$ denoting  an estimator  of  $R$, the prediction error ${\mathcal E}_{n+h}$ can be decomposed as 
\begin{equation*}
 \X_{n+h}-\widehat{\X}_{n+h}  = \varepsilon_{n+h} + \big(\Phi_1(\X_{n+h-1})-R(\X_{n+h-1}) + \Phi_2(\X_{n+h-2})\big) + \big(R(\X_{n+h-1})-\widehat{R}(\widehat{\X}_{n+h-1})\big).
 \end{equation*}
Consider now the conditional distribution  ${\mathcal E}_{n+h} |\X_n$. Notice that if $ h=1$, the model specification error $\big(\Phi_1(\X_n)-R(\X_{n}) + \Phi_2(\X_{n-1})\big)$ causes a shift in this conditional distribution  due to  the term $\Phi_1(\X_n)-R(\X_{n})$ as well as an increase in  variability due to the term $ \Phi_2(\X_{n-1})$. Similarly,  for $h \geq 2$ the model specification error $\big(\Phi_1(\X_{n+h-1})-R(\X_{n+h-1}) + \Phi_2(\X_{n+h-2})\big)$ does not vanish asymptotically.
Furthermore,    for $h\geq 2$,  the  error $\big(R(\X_{n+h-1})-\widehat{R}(\widehat{\X}_{n+h-1})\big)$ is not only due to the estimator $ \widehat{R}$ of $ R$ (as in the case $h=1$), but also due the fact that the unknown random element $ \X_{n+h-1}$ has been replaced  by  its predictor $ \widehat{\X}_{n+h-1}$.  This causes a further increase in variability. 

An appropriate procedure to construct prediction intervals or bands, should consider {\it all three} aforementioned sources affecting the prediction error and consistently estimate the conditional distribution $\mathcal E_{n+h} | \X_n , \X_{n-1}, \ldots, \X_{n-k+1}$. However, and to the best of our knowledge, this issue has not been appropriately explored in the literature. In particular, and even in the most studied univariate, real-valued case, it is common to estimate the prediction error distribution by ignoring the model specification error, that is,  by assuming that the model used for prediction is identical to the data generating process. Consequently, the bootstrap approaches applied in this context, use the same model to make the prediction and to generate the bootstrap pseudo-time series. Such approaches ignore the model misspecification error; see \cite{TS90}, \cite{BDD95}, \cite{APR02}, \cite{PRR04} as well as \cite{PP16b} and the references therein. See also Section~\ref{sec:3} for more details.

In this paper, we develop a bootstrap procedure to construct prediction bands for functional time series that appropriately takes into account  {\it all three} sources of errors affecting the conditional distribution of ${\mathcal E}_{n+h}$. The proposed bootstrap approach generates, in a model-free way, pseudo-replicates $\X_1,^*, \X_{2}^*, \ldots, \X_{n}^*$, and $ \X_{n+1}^* , \X_{n+2}^*, \ldots, \X_{n+h}^*$ of the functional time series at hand that  appropriatelly mimic the dependence structure of the underlying functional process. Moreover, the approach ensures that the generated functional pseudo-time series has the same $k$ functions at the end as the functional times series observed; that is, $ \X_t^*=\X_t$ holds true for $ t\in\{n-k+1, n-k+2, \ldots, n\}$. This is important because, as  already mentioned, it is the conditional distribution of $\mathcal E_{n+h} $ given $ \X_n, \X_{n-1}, \ldots, \X_{n-k+1}$ in which we are interested. These  requirements are fulfilled by generating  the functional pseudo-elements  $ \X_1^*, \X_{2}^*, \ldots, \X_{n}^*$ using a backward-in-time vector autoregressive representation of the  time-reversed  process of scores appearing in the  Karhunen-Lo\`{e}ve representation (see Section~\ref{sec:2} for details). Given the model-free, bootstrap-generated functional pseudo-time series   $ \X_1^*, \X_{2}^*, \ldots, \X_{n}^*$ and $ \X^*_{n+1}, \ldots, \X_{n+h}^*$, the same  model used to obtain  the predictor $ \widehat{\X}_{n+h}=\widehat{g}(\widehat{\X}_{n+h-1}, \ldots, \widehat{\X}_{n+h-k})$, see~\eqref{eq.pred1},  is then applied,  and the  pseudo-predictor  $\widehat{\X}^*_{n+h} =\widehat{g}^*(\widehat{\X}^*_{n+h-1}, \ldots, \widehat{\X}^*_{n+h-k})$ is obtained. Here,  $\widehat{\X}^*_t =\X^*_t=\X_t$ if $ t \in \{n, n-1, \ldots, n-k+1\}$ and  $ \widehat{g}^*$ denotes the same estimator as $\widehat{g}$ but based on the  bootstrap functional pseudo-time series $\X_1^*, \X_{2}^*, \ldots, \X_{n}^* $. The conditional (on ${\mathcal X}_{n,k}$) distribution of the prediction error $\X_{n+h}- \widehat{\X}_{n+h}$ is then estimated using the conditional distribution of the bootstrap prediction error $\X^*_{n+h}- \widehat{\X}^*_{n+h}$. We show that the described procedure leads to consistent estimates of the conditional distribution of interest. We also prove the consistency of the bootstrap in estimating the conditional distribution of the studentized prediction error process in $ {\mathcal H}$. The latter consistency is important because it  theoretically justifies the use of the proposed bootstrap method in the construction of simultaneous prediction bands for the $h$-step-ahead prediction  that also appropriately account for the local variability of the corresponding prediction error. Using simulations and two empirical data applications the good finite sample performance of the  bootstrap procedure proposed.

We perform dimension reduction via the truncation of the Karhunen-Lo\`{e}ve representation and we capture the infinite-dimensional structure of the underlying functional process  by allowing the number of principal components used to increase to infinite (at an appropriate rate) with $n$. These aspects of our procedure  are common to the sieve-bootstrap introduced in \cite{Paparoditis18}. However, apart from the differences in the technical tools used for establishing bootstrap validity, a novel, general backward autoregressive representation of the vector process of scores is introduced, which is a key part of the bootstrap procedure proposed in this paper. This representation allows for the generation of the functional pseudo time series $\X_1^\ast, \X_2^\ast, \ldots, \X_n^\ast$, backward in time and enables, therefore,  for this pseudo time series  to satisfy the condition $\X_t^*=\X_t$ for $t\in \{n-k+1, n-k+2, \ldots, n\}$. The latter condition is essential for successfully evaluating the conditional distribution of the prediction error $ {\mathcal E}_{n+h} |{\mathcal X}_{n,k}$ using the bootstrap and for the construction of the desired prediction bands. The same condition, including the aforementioned backward vector autoregressive representation of the score process and the focus on prediction error distribution, also are the main differences to the resampling approach considered in \cite{Shang18} which has been used for estimating the long-run variance. This approach is based on bootstrapping the principal component scores by maximum entropy.  

 \cite{APS06} and \cite{ABC+16} considered nonparametric (kernel-type), one step ahead predictor and proposed a resampling method to construct pointwise prediction intervals. Model-based bootstrap procedures to construct pointwise prediction intervals for one step ahead prediction under (mainly) FAR(1) model assumptions on the data generating process have been considered in \cite{RAV+16} and \cite{VAR18}. These approaches differ from ours. They are developed for particular predictors, they are designed for pointwise prediction intervals only, and the bootstrap approaches involved, only are  designed   for the specific prediction setting  considered. More related to our approach in terms of not requiring particular model assumptions for the data generating process and of not being designed for a specific predictor,  is the approach proposed in \cite{ANH15} for the construction of prediction bands. In Section~\ref{sec:5} we compare the performance of this approach with the bootstrap approach proposed in this paper. 

This paper is organized as follows. In Section~\ref{sec:2}, we state the notation used  and introduce the notion of backward vector autoregressive representations of the time-reversed vector process of scores appearing in the Karhunen-Lo\`{e}ve representation. In Section~\ref{sec:3}, we present the proposed bootstrap procedure and show its asymptotic validity for the construction of  simultaneous prediction bands. Section~\ref{sec:4} is devoted to some practical issues related to the construction of prediction bands and the implementation of our procedure. Section~\ref{sec:5} investigates the finite sample performance of the proposed bootstrap procedure using simulations, while in Section~\ref{sec:6}, applications of the new methodology to a real-life data set are considered. Conclusions are provided in Section~\ref{sec:7}. The analysis of a second real-life data set as well 
as the proofs of the theoretical results presented in this paper together with some auxiliary lemmas are given  in the Supplementary Material.

\section{Preliminaries}\label{sec:2}

\subsection{Setup and Examples of Predictors}\label{sec:2.1}

Consider a time series  $\X_1, \X_2, \dots, \X_n$  stemming from a stationary,   $ L^4$-${\mathcal M}$-approximable stochastic process ${\bf X}=\{\X_t, t \in \mathbb{Z}\}$ with mean $E (\X_t)=0$  and autocovariance operator $C_r = E(\X_t\otimes \X_{t+r})$,  $r\in\mathbb{Z}$. Recall that $C_r$ is a Hilbert-Schmidt (HS) operator.  The  $L^4$-${\mathcal M}$ approximability property allows for a weak dependence structure of the underlying functional process, which covers a wide range of commonly used functional time series models, including functional linear processes and functional autoregressive,  conditional heteroscedasticity processes; \citep[see][for details]{HK10}. $L^4$-${\mathcal M}$-approximability   implies that $\sum_{r\in Z}\|C_r\|_{\text{HS}} <\infty$ and therefore, that the functional process $ {\bf X}$ possesses  a continuous, self-adjoint spectral density operator ${\mathcal F}_\omega$,   given by 
\begin{equation*}
{\mathcal F}_\omega =(2\pi)^{-1} \sum_{r\in Z} C_h e^{-ir \omega}, \omega \in \mathbb{R},  
\end{equation*}
which is trace class \citep{HKH15} (also see \cite{PT13} for a different set of weak dependence conditions on the functional process ${\bf X}$). Here and in the sequel, $ \|\cdot\|_{HS}$ denotes the Hilbert-Schmidt norm of an operator while $ \|\cdot \|_F$ the Frobenious norm of a matrix. We assume that the eigenvalues  $ \nu_1(\omega), \nu_2(\omega), \ldots, \nu_m(\omega)$ of the spectral density operator $ {\mathcal F}_\omega$ are strictly positive for every $\omega \in [0,\pi]$. 

Suppose that the h-step-ahead predictor of $ \X_{n+h}$ is  obtained as 
\begin{equation} \label{eq.xpred}
\widehat{\X}_{n+h} = \widehat{g}_{(h)}(\X_{n},\ldots, \X_{n-k+1}),
\end{equation}
where $ k \in N$, $k<n$,  is fixed and determined by the model  selected to perform the  prediction (see  also~\eqref{eq:model}),   while  $ \widehat{g}_{(h)}$ denotes  an estimator of the unknown operator $g_{(h)}$. Based on  $\widehat{\X}_{n+h}$,  our aim is to construct a prediction interval, respectively,  prediction band  for $ \X_{n+h}$ associated with the model~\eqref{eq.xpred} which is used for prediction. Toward this end,   an  estimator of the  distribution of the prediction error $\mathcal E_{n+h} = \X_{n+h}-\widehat{\X}_{n+h}$ is needed. More precisely, we are interested in estimating the conditional distribution 
\begin{align} \label{eq.prederr}
\mathcal E_{n+h}  \big| \X_n, \X_{n-1}, \ldots, \X_{n-k+1}. 
\end{align}
Since  we do not want to restrict our considerations to a specific predictor $\widehat{g}$, many of the  predictors applied  in the functional time series literature  fit in our setup. We elaborate on  some  examples:
\begin{enumerate}
\item[1)] Suppose that  in~\eqref{eq:model}  the operator $g$ is given by $ g(\X_{n},\ldots, \X_{n-k+1}) =\sum_{j=1}^k \Phi_j(\X_{n+1-j})$ with the $ \Phi_j $'s  being  linear, bounded operators  $ \Phi_j :{\mathcal H} \rightarrow {\mathcal H}$. This is a case where  a functional autoregressive model of order $k$ (FAR($k$)) is used to predict $ \X_{n+h}$, see \citep{KR13} in which the issue of  the selection of the  order k  also  is  discussed. Given  some estimators $ \widehat{\Phi}_j$ of $ \Phi_j$,  the corresponding $h$-step-ahead predictor  is given by $  \widehat{g}(\widehat{\X}_{n+h-1},\ldots, \widehat{\X}_{n+h-k})=\sum_{j=1}^k\widehat{\Phi}_j(\widehat{\X}_{n+h-j})$, where $ \widehat{\X}_t\equiv \X_t$ if $ t \in \{n, n-1, \ldots, n-k+1\}$.  A special case  is  the popular FAR(1)  model in which  it is assumed that   $\X_t$  is generated as  $\X_t=\Phi(\X_{t-1}) + v_t$ with  $ \|\Phi\|{_\mathcal L} <1$ and $ \varepsilon_t $ an i.i.d. sequence in $ {\mathcal H}$ \citep{Bosq00, BB07}. Here  and in the sequel, $ \|\cdot\|_{\mathcal L}$ denotes the operator norm. 
\item[2)]  Suppose that  $ g_{(h)}(\X_n,\ldots, \X_{n-k+1}) = \sum_{j=1}^ d {\bf 1}^{\top}_j\sum_{l=1}^kD_l\bm\xi_{n+h-l} v_j$,  where ${\bf 1}_j $ is the $d$-dimensional vector with the $j$\textsuperscript{th} component equal to 1 and 0 elsewhere,  $ \bm\xi_t$ is the $d$-dimensional vector $\bm\xi_t= (\langle \X_t,v_j\rangle, j=1,2, \ldots, d)^{\top}$, $v_j$ are the orthonormal eigenfunctions corresponding to the $d$ largest eigenvalues of the lag-0 covariance operator $C_0=E(\X_0 \otimes \X_0)$, and ($ D_1, D_2, \ldots, D_k$) are the matrices obtained by the orthogonal projection of $ {\bm \xi}_t$ on the space spanned by ($\bm\xi_{t-1}, \bm\xi_{t-2}, \ldots, \bm\xi_{t-k}$). A predictor of $ \X_{n+h}$ is then obtained as 
\begin{equation*}
 \widehat{\X}_{n+h}= \sum_{j=1}^d{\bf 1}^{\top}_j \breve{\bm\xi}_{n+h}\widehat{v}_j, 
\end{equation*}
where $\breve{\bm\xi}_{n+h}=\sum_{l=1}^k \widehat{D}_l\breve{\bm\xi}_{n+h-l}$,  $\breve{\bm\xi}_{t}=  \widehat{\bm\xi}_t$ for  $ t \in \{n,n-1, \ldots, n-k+1\}$ and  $ \widehat{\bm\xi}_1,  \ldots, \widehat{\bm\xi}_n$ are  the estimated $d$-dimensional score vectors
 $ \widehat{\bm\xi}_t = (\langle \X_t,\widehat{v}_j\rangle , j=1,2, \ldots, d)^{\top}$. Here $ \widehat{v}_j$ are the estimated orthonormal eigenfunctions  corresponding to the $d$ largest estimated  eigenvalues  of $ \widehat{\mathcal C}_0=n^{-1}\sum_{t=1}^n (\X_t-\overline{\X}_n)\otimes  (\X_t-\overline{\X}_n)$ and ($ \widehat{D}_l$, $l=1,2, \ldots, k$) are the estimated $d\times d$ matrices obtained by least squares fitting of a $k$\textsuperscript{th} order vector autoregression to the time series  $ \widehat{\bm\xi}_t$, $ t=1,2, \ldots, n$ \citep{ANH15}. 
\item[3)] Similar to 2), the predictor of $\X_{n+h}$  can be  obtained as $\widehat{\X}_{n+h}= \sum^d_{j=1}{\bf 1}_j^{\top}\breve{\bm{\varepsilon}}_{n+h,j}\widehat{v}_j$, where $\bm{\breve{\varepsilon}}_{n+h,j}$ is  an $h$-step predictor of the $j$\textsuperscript{th} component,  obtained via a univariate time series forecasting method applied to estimated components  $(\widehat{\varepsilon}_{1,j}, \dots, \widehat{\varepsilon}_{n,j})$ for each $j=1,\dots,d$, \citep{HS09}.  
 \item[4)] Let for notational simplicity $k=1$ and  $g_{(h)}(\X_n)=E(\X_{n+h}|\X_n)$  be the conditional mean function of $ \X_{n+h}$ given $ \X_n$. Consider the  predictor $\widehat{\X}_{n+h}$  obtained  using a nonparametric estimator of $g_{(h)}$,   for instance, a functional version of the Nadaraya-Watson estimator given by 
\begin{align*}
\widehat{g}_{(h)}(\X) &= 
 \sum^{n-h}_{i=1}\frac{K\left[d(\X_i, \X)/\delta\right]\X_{i+h}}{\sum^{n-1}_{j=1}K\left[d(\X_j, \X)/\delta\right]},
\end{align*}
where $K(\cdot)$ is a kernel function,  $\delta>0$ is a smoothing bandwidth,  and $d(\cdot, \cdot)$ is a distance function on $ {\mathcal H}$. $ \widehat{\X}_{n+h}=\widehat{g}_{(h)}(\X_n)$ is then the predictor of $\X_{n+h}$, \citep[see, e.g.,][]{APS06}.
\end{enumerate}

\subsection{The Time-Reversed Process of Scores}

To introduce the proposed bootstrap procedure, it is important to first discuss some properties of the time-reversed process of scores associated with the functional process ${\bf X}$.  To this end, consider for  $m\in \mathbb{N}$,   the   $m$-dimensional vector process of scores, that is   $ \bm\xi=\{ \bm\xi_t, t \in \mathbb{Z}\}$,  where $ \bm\xi_t=(\xi_{j,t}=\langle \X_t,v_j \rangle, j=1,2, \ldots, m)^{\top}$ and  $ v_1, v_2, \ldots $, are the orthonormal eigenvectors corresponding to the eigenvalues $\lambda_1>\lambda_2 > \ldots$,  in descending  order, of the lag-0 autocovariance operator $ C_0$.  Denote by $ \widetilde{\bm\xi}= \{\widetilde{ \bm\xi}_{t}, t \in \mathbb{Z}\}$ the time-reversed version of $  \bm\xi $, that is, $\widetilde{\bm\xi}_{t}= \bm\xi_{-t}$ for any  $t \in \mathbb{Z}$. We call $\bm\xi$ and $\widetilde{\bm\xi}$ the forward and the backward score processes, respectively. The autocovariance structure of both processes is closely related because for any $h\in \mathbb{Z}$ we have
\begin{align} \label{eq.autoc}
\Gamma_{\bm{\widetilde{\xi}}}(h) & := E\left[\bm{\widetilde{\xi}}_0(m)\bm{\widetilde{\xi}}_h^{\top}(m)\right]  \nonumber \\
& = E\left[\bm{\xi}_{0}(m)\bm{\xi}_{-h}^{\top}(m)\right] =:  \Gamma_{\bm{\xi}}(-h).
\end{align}
Thus, properties of the forward score process $\bm\xi$,  which arise from its second-order structure,  carry over to the backward process $ \widetilde{\bm\xi}$.  To elaborate, note first that  the (Hilbert-Schmidt) norm summability of the autocovariance operators $C_h$  as well as the assumption that  the eigenvalues  $ \nu_1(\omega), \nu_2(\omega), \ldots, \nu_m(\omega)$ of the spectral density operator $ {\mathcal F}_\omega$ are  bounded away from zero for all $\omega \in  [0,\pi]$,  imply that,  the $m\times m$ spectral density matrix $ f_{\bm\xi}(\omega) =(2\pi)^{-1}\sum_{h\in Z}\Gamma_{\bm\xi}(h)e^{-ih\omega}$  of the forward score process  $ \bm\xi$, is  continuous,  bounded from above and  bounded  away from zero from below; 
\citep[see Lemma 2.1 of][]{Paparoditis18}. The same properties also hold true for the      $m\times m$ spectral density matrix 
$ f_{\widetilde{\bm\xi}}(\omega)=(2\pi)^{-1}\sum_{h\in Z}\Gamma_{\widetilde{\bm\xi}}(h)e^{-ih\omega}$ of the backward score process $ \widetilde{\bm\xi}$.  This follows  immediately from the corresponding  and aforementioned  properties of  $f_{\bm\xi} $  and taking  into account that by equation \eqref{eq.autoc},  $f_{\widetilde{\bm\xi}}(\omega)=f^{\top}_{\bm\xi}(\omega)$ for all $\omega \in [0,\pi]$. Now, the fact that both spectral density matrices $ f_{{\bm\xi}}$ and  $f_{\widetilde{\bm\xi}}$ are bounded from above and from bellow, implies by Lemma 3.5 of \citet[][p. 116]{CP93},  that both processes--the process  ${\bm\xi} $ and the time-reversed process $\widetilde{\bm\xi} $--obey a so-called vector autoregressive representation. That is,  infinite sequences of $m\times m$ matrices $\{ A_j, j \in N\}$ and $ \{B_j, j \in N\}$ as well as  full rank $m$-dimensional, white noise processes $ \{\bm{e}_t, t \in Z\}$ and $ \{\bm{v}_t, t \in Z\} $  exist such that the random vectors  $\bm\xi_t$  and $ \widetilde{\bm\xi}_t$ have, respectively, the following  autoregressive representations:    
\begin{equation}\label{eq.var-forw}
\bm\xi_t = \sum_{j=1}^\infty A_j \bm\xi_{t-j} + \bm{e}_t
\end{equation}
and
\begin{equation}\label{eq.var-back}
\tilde{\bm{\xi}}_t = \sum^{\infty}_{j=1}B_j\tilde{\bm{\xi}}_{t-j} + \bm{u}_t.
\end{equation}
We refer to~\eqref{eq.var-forw} and~\eqref{eq.var-back} as the forward and backward vector autoregressive representations of ${\bm \xi}_t$ and to $ \{\bm{e}_t\}$ and to $\{\bm{u}_t\}$ as the  forward and the backward noise processes. We stress here the fact that representations~\eqref{eq.var-forw} and~\eqref{eq.var-back} should  not be confused with that of a {\it linear} (infinite order)  vector autoregressive process. This is due to the fact that the noise vector processes   $ \{\bm{e}_t\}$ and  $\{\bm{u}_t\}$  appearing in representations~\eqref{eq.var-forw} and~\eqref{eq.var-back}, respectively,  are only uncorrelated and not necessarily i.i.d. sequences of random vectors. Furthermore, the autoregressive  matrices $\{A_j\}$ and $ \{B_j\}$ appearing in the above representations also satisfy the summability conditions $ \sum_{j=1}^{\infty} \|A_j\|_F <\infty$ and $  \sum_{j=1}^{\infty} \|B_j\|_F <\infty$, while the corresponding power series 
\begin{align*}
A(z)=I-\sum_{j=1}^\infty A_j z^{j} \ \ \ \mbox{and} \ \ \ 
B(z)=I-\sum_{j=1}^{\infty} B_j z^{j}
\end{align*}
do not vanish for $ |z| \leq 1$; that is, $ A^{-1}(z)$ and $ B^{-1}(z)$ exist for all $ |z| \leq 1$ \citep[see][for more details on such vector autoregressive representations of weakly stationary processes]{CP93,MK15}. Using reversion in time and, specifically, the property that $\widetilde{\bm\xi}_t=\bm\xi_{-t}$, equation~\eqref{eq.var-back} leads to the expression 
\begin{equation} \label{eq.revers}
 \bm{\xi}_{t}=\sum^{\infty}_{j=1}B_j\bm{\xi}_{t+j}+\bm{u}_{t}, 
\end{equation}
which also can be written as $ B(L^{-1})\bm{\xi}_t = \bm{u}_t$, with the shift operator $L$ defined by  $ L^k \bm{\xi}_t = \bm{\xi}_{t-k}$ for any  $ k \in \mathbb{Z}$. Expression~\eqref{eq.revers} implies that the two white noise  innovation processes $\{\bm{e}_t,t\in \mathbb{Z}\}$ and $ \{\bm{u}_t, t\in \mathbb{Z}\}$, are related by 
\begin{align} \label{eq.v-e}
\bm{u}_t  & = B(L^{-1})\bm\xi_t =B(L^{-1})A^{-1}(L) \bm{e}_t, \ \ \  t\in Z.
\end{align} 
Notice that~\eqref{eq.v-e} generalizes to the vector autoregressive case an analogue  expression obtained for the univariate autoregressive case by   \cite{Findley86} and \cite{BDD95}. Further, and as relation~\eqref{eq.v-e} verifies, even if $\bm\xi_t$ in~\eqref{eq.var-forw} is a linear process, that is  even if $\{ \bm{e}_t\}$  is an i.i.d. innovation process in $\mathbb{R}^m$, the  white noise innovation process $ \{\bm{u}_t\}$ appearing in the  time-reversed process~\eqref{eq.revers} is, in general, not an i.i.d. process.  

\section{Bootstrap Prediction Intervals}\label{sec:3}

\subsection{Bootstrap Procedure}

The basic idea of the proposed bootstrap procedure is to generate a functional time series of pseudo-random elements $\X_1^*, \X_2^*, \dots, \X_n^*$, and future values $ \X_{n+1}^*, \X_{n+2}^*, \ldots, \X_{n+h}^*$, which appropriately imitate the dependence structure of the functional time series at hand, while at the same time satisfy the condition
\begin{equation} \label{eq.requir}
  \X^\ast_{n-k+1}=\X_{n-k+1}, \quad  \X^\ast_{n-k+2}=\X_{n-k+2},\  \ldots,  \  \X^\ast_{n}=\X_{n}.
\end{equation}
The above condition is important because, as we have seen, the conditional distribution of $\mathcal E_{n+1}(\cdot)$ given  $\X_{n}, \X_{n-1}, \ldots, \X_{n-k+1}$ is the one in which we are interested. Toward this goal and motivated by the functional sieve bootstrap proposed by \cite{Paparoditis18}, we use the Karhunen-Lo\`{e}ve representation  and decompose  the random element $ \X_t$ in two parts:
\begin{equation} \label{eq.dec}
\X_t = \sum^{\infty}_{j=1}\xi_{j,t}v_j =\underbrace{\sum^m_{j=1}\xi_{j,t}v_j}_{\displaystyle \X_{t,m}} + \underbrace{\sum_{j=m+1}^{\infty}\xi_{j,t}v_j}_{\displaystyle U_{t,m}}.
\end{equation} 
In \eqref{eq.dec}, the element $\X_{t,m}$ is considered as the main  driving part of $ \X_t$, while    the ``remainder'' $U_{t,m}$ is treated as  a white noise component. Now, to generate  the functional pseudo-time series $\X_1^*, \X_2^*, \dots, \X_n^*$, we first bootstrap the $m$-dimensional time series of scores  by  using the backward vector autoregressive representation  given in~\eqref{eq.revers}. Using the backward  representation  allows for the generation of a  pseudo-time series of scores $ \bm\xi_1^*,  \bm\xi_2^*, \ldots,  \bm\xi_n^* $, which  satisfies the condition $ \bm\xi_{t}^*=\bm\xi_t$ for $ t\in\{n-k+1, n-k+2, \ldots, n\}$. This is important to ensure  that  the bootstrap-generated time series $ \X_1^*, \X^*_2, \ldots, \X_n^*$  fulfills requirement~\eqref{eq.requir}. The backwards-in-time-generated pseudo-time series of  scores  $ \bm\xi_1^*, \bm\xi_2^*, \ldots, \bm\xi_n^*$  can then be transformed to pseudo-replicates of the main driving part $ \X_{t,m}$ by using the equation   $ \X_{t,m}^\ast=\sum^m_{j=1}\xi_{j,t}^\ast v_j$. Notice that since by construction  $\bm\xi^\ast_t=\bm\xi_t$ for  $ t=n, n-1, \ldots, n-k+1$,  we  have that  $\X^*_{t,m}=\X_{t,m} $ and, consequently,  we set  $  \X^\ast_t=\X_t$ for the same set of time indices.  Adding to the generated $ \X_{t,m}^\ast$  for  the remaining indices $ t=n-k,n-k-1,  \ldots, 1$, an appropriately resampled functional noise $ U^\ast_{t,m}$,  leads to the functional pseudo replicates $ \X_1^*, \X_2^*, \dots, \X_{n-k}^*$. As a result, a  functional pseudo-time series $ \X_1^*, \X_2^*, \ldots, \X_n^* $ can be obtained that imitates the dependence structure of $ \X_1, \X_2, \ldots, \X_n$ and  at the same time satisfies~\eqref{eq.requir}. Notice that  implementation of the above ideas  requires estimation of the  eigenvectors $v_j$ and of the scores $\xi_{j,t}=\langle \X_t, v_j\rangle$,  because these quantities are not observed (see Section~\ref{sec.bootalgo} for details). 

\begin{figure}[tbp]
\centering
{\includegraphics[width=9cm]{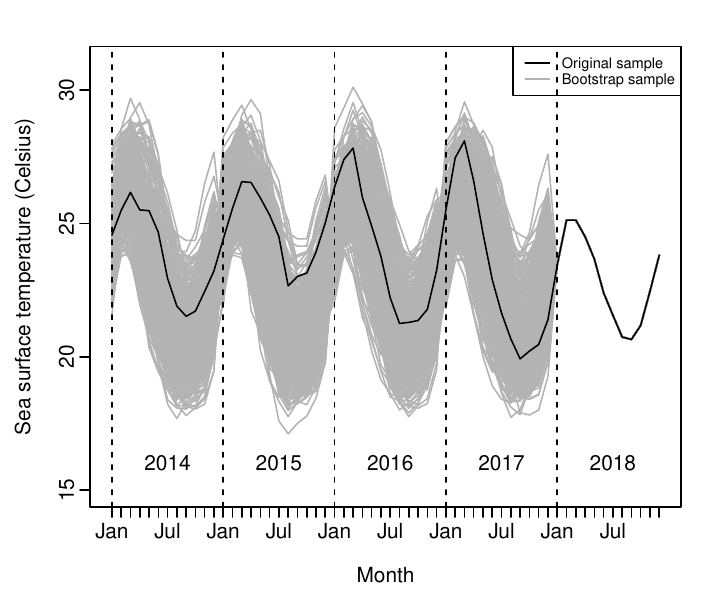}}
\caption{Sea surface temperature in El Ni\~{n}o region 1+2 displayed from January 2014 to December 2018 (black line) together with 1000 different bootstrap samples (gray lines) when $k=1$.}\label{fig:0}
\end{figure}

Before proceeding with a precise description of the bootstrap algorithm, we illustrate its capability using a data example. Figure~\ref{fig:0} shows the monthly sea surface temperatures for the last three years (analyzed in Section~\ref{sec:6}) together with 1,000 bootstrap replications obtained when $k=1$ and using the bootstrap algorithm described in Section~\ref{sec.bootalgo}. Notice the asymmetric features of the time series paths generated and the fact that all 1,000 bootstrap samples displayed pass through the same final curve. That is, all generated bootstrap functional time series satisfy condition~\eqref{eq.requir}, which for the case $k=1$ reduces to $ \X^*_n=\X_n$.

\subsection{Bootstrap Algorithm}

We now proceed with a detailed description of the bootstrap algorithm used  to generate the functional pseudo-time series $ \X_1^*, \X_2^*, \ldots, \X_n^*$ and $ \X_{n+1}^\ast, \X_{n+2}^*, \ldots, \X^*_{n+h}$. Steps~1 to~3 of the following algorithm concern the generation of the future pseudo-elements $ \X^\ast_{n+1},\X_{n+2}^*,  \ldots,  \X^*_{n+h}$, while Steps~4 to~6  the generation of $ \X_1^\ast, \ldots, \X^\ast_{n}$.

\vspace*{0.2cm}

\label{sec.bootalgo}
\begin{enumerate}
\item[] {\bf Step~1:} Center the observed functional time series by calculating  $ {\mathcal Y}_t={\mathcal X}_t-\overline{\mathcal X}_n$, $ \overline{\mathcal X}_n=n^{-1}\sum_{t=1}^n{\mathcal X}_t$.
\item[] {\bf Step~2:} Select integers  $m $ and $p$, where $m$ is the   truncation number in~\eqref{eq.dec} and $ p$  the  order used to approximate the infinite-order vector autoregressive representations~\eqref{eq.var-forw} and~\eqref{eq.var-back}. Denote by $\widehat{\bm\xi}_1, \widehat{\bm\xi}_2, \ldots, \widehat{\bm\xi}_n $ the time series of estimated,  $m$-dimensional vector of scores; that is, 
\[  \widehat{\bm\xi}_t =\big(\langle {\mathcal Y}_t,\widehat{v}_j\rangle, \ j=1,2, \ldots, m\big)^\top, \ t=1,2, \ldots, n \]
 where $ \widehat{v}_j$, $ j=1,2, \ldots, m$,   are the  estimated (up to a sign) orthonormal eigenfunctions corresponding to the $m$ largest 
estimated eigenvalues of the lag-0 sample autocovariance operator $ \widehat{C}_0=n^{-1}\sum_{t=1}^n {\mathcal Y}_t \otimes {\mathcal Y}_t$.
\item[] {\bf Step~3:}  Fit a VAR$(p)$ process to the ``forward'' series of estimated scores; that is,
$\widehat{\bm{\xi}}_t = \sum^p_{j=1}\widehat{A}_{j,p}\widehat{\bm{\xi}}_{t-j}+\widehat{\bm{e}}_t$,  $
 t=p+1,p+2,\dots, n$,
with $ \widehat{\bm{e}}_t$ being the estimated residuals.  Generate   
$ \bm{\xi}^*_{n+h}= \sum_{l=1}^{p} \widehat{A}_{l,p}{\bm{\xi}}^\ast_{n+h-l} + \bm{e}^*_{n+h}$,
where we set ${\bm{\xi}}^\ast_{n+h-l}=\widehat{\bm{\xi}}_{n+h-l} $  if $n+h-l\leq n $ and $ \bm{e}^*_{n+h} $ is i.i.d. resampled from the set of centered residuals $ \{\widehat{\bm{e}}_t-\overline{\bm{e}}_n , t=p+1, p+2, \ldots, n\}$, $ \overline{{\bm e}}_n = (n-p)^{-1}\sum_{t=p+1}^n \widehat{{\bm e}}_t$. 
Calculate 
\begin{align*}
 \X^*_{n+h} &= \overline{\mathcal X}_n +  \sum_{j=1}^m{\bf 1}_{j}^\top \bm{\xi}^*_{n+h} \widehat{v}_j +  U^*_{n+h,m},
 \end{align*}
where the $ U^*_{n+h,m}$ are   i.i.d. resampled from the set   $\big\{\widehat{U}_{t,m} - \overline{U}_n,$ $t=1,2,\dots,n\big\}$,  $\overline{U}_n = n^{-1}\sum^n_{t=1}\widehat{U}_{t,m}$ and $\widehat{U}_{t,m} = {\mathcal Y}_t - \sum^m_{j=1}{\bf 1}_j^{\top}\widehat{\bm\xi}_{t}\widehat{v}_j$. Recall that  ${\bf 1}_j $  denotes the $m$-dimensional vector with the $j$\textsuperscript{th} component equal to 1 and 0 elsewhere. 

If $ p \leq k+h$,  move to Step~4. If $ p >k+h$,  generate  for $l=1,2, \ldots, p-(k+h)$  additional random vectors $ \bm\xi^*_{n+h+l}=\sum_{j=1}^{p} \widehat{A}_{j,p}\bm\xi^{*}_{n+h+l-j} + \bm{e}^\ast_{n+h+l}$, where $ \bm{e}^*_{n+h+l}$ are i.i.d. generated in the same way as $  \bm{e}^*_{n+h}$.
\item[] {\bf Step~4:} Fit a VAR$(p)$ process to the ``backward" series of estimated scores; that is,
\begin{equation*}
\widehat{\bm{\xi}}_t= \sum^p_{j=1}\widehat{B}_{j,p}\widehat{\bm{\xi}}_{t+j}+\widehat{\bm{u}}_t, \qquad t=1,2,\dots,n-p.
\end{equation*}
\item[] {\bf Step~5:} Generate a pseudo-time series of the scores $\{\bm{\xi}_1^*, \bm{\xi}_2^*,\dots, \bm{\xi}_n^*\}$ by setting
$\bm{\xi}_t^* = \widehat{\bm\xi}_t$ for 
$ t=n, n-1,\dots, n-k+1$,	
and  by using  for $ t=n-k, n-k-1,\dots,1$, the backward vector autoregression 
$\bm{\xi}_t^* = \sum^p_{j=1}\widehat{B}_{j,p}\bm{\xi}_{t+j}^* + \bm{u}^*_t$.
Here $\bm{u}_1^*, \bm{u}_2^*,\ldots,\bm{u}_{n-k}^* $ are obtained as (see~\eqref{eq.v-e}), 
\[ \bm{u}^*_t = \widehat{B}_p(L^{-1})\widehat{A}_p^{-1}(L) \bm{e}^*_t,\]
with  $\widehat{A}_p(z)=I-\sum_{j=1}^p \widehat{A}_{j,p}z^j$, $ \widehat{B}_p(z)=I-\sum_{j=1}^p \widehat{B}_{j,p}z^{j}$, $z\in {\mathbb C}$, and where the $\bm{e}_t^*$ are i.i.d. resampled as in Step~3. 
\item[] {\bf Step~6:} Generate a pseudo-functional time series $\{\X_1^*, \X_2^*, \dots, \X_n^*\}$ as follows. For $  t=n, n-1, \dots, n-k+1$ set  
\begin{align*}
\X_t^* &=\overline{\mathcal X}_n +  \sum_{j=1}^m{\bf 1}_j^\top \widehat{\bm{\xi}}_t\widehat{v}_j + \widehat{U}_{t,m}  \  \equiv \ \X_t, 
\end{align*}
while for $ t=n-k, n-k-1, \ldots, 1$, use the obtained backward pseudo-scores $ \bm\xi_1^*, \bm\xi_2^*, \ldots, \bm\xi_{n-k}^*$ and  calculate 
\begin{equation*} \label{eq.xstar}
\X_t^* = \overline{\mathcal X}_n + \sum^m_{j=1}{\bf 1}^{\top}_j\bm{\xi}_t^*\widehat{v}_j + U_{t,m}^*.
\end{equation*}
Here, the $U_{t,m}^*$ are i.i.d.  pseudo-elements resampled as in Step~3. 
\item[] {\bf Step~7:} If   model~\eqref{eq:model} is used to obtain  the prediction $\widehat{\X}_{n+h}$, then  calculate the pseudo-predictor
\begin{equation} \label{eq.pred-boot}
\widehat{\X}_{n+h}^* =  \overline{\mathcal X}_n + \widehat{g}^\ast\left( \widehat{\X}^*_{n+h-1}-\overline{\X}^*_n, \widehat{\X}^\ast_{n+h-2}-\overline{\X}^*_n,  \dots, \widehat{\X}^*_{n+h-k}- \overline{\mathcal X}^*_n\right),
\end{equation}
where we set $ \widehat{\X}^*_{t}-\overline{\X}_n^\ast \equiv \X_t-\overline{\X}_n$ for  $ t=n,n-1, \ldots, n-k+1$, $ \overline{X}_n^*=n^{-1}\sum_{t=1}^n\X^*_t$  and $\widehat{g}^*$ is the same estimator as $\widehat{g}$ used in~\eqref{eq.pred1} but obtained using the  generated pseudo-time series $ \X_1^*, \X_2^*, \dots, \X_n^*$. Alternatively, the bootstrap analogue of~\eqref{eq.xpred}  can be  calculated as
\begin{equation} \label{eq.pred-boot2}
\widehat{\X}_{n+h}^* =  \overline{\mathcal X}_n + \widehat{g}^\ast_{(h)}\left( \X_n-\overline{\X}_n, \X_{n-1}-\overline{\X}_n,  \dots, \X_{n-k+1}- \overline{\mathcal X}_n\right).
\end{equation}
Here $\widehat{g}^*_{(h)}$ is the same estimator as $\widehat{g}_{(h)}$  given in~\eqref{eq.xpred} but  based on the pseudo-time series $ \X_1^*, \X_2^*, \dots, \X_n^*$.
\item[] {\bf Step~8:}  Use the distribution of  $\mathcal E^\ast_{n+h}=  \X^\ast_{n+h}-\widehat{\X}^\ast_{n+h}$ to approximate  the (conditional) distribution of $\mathcal E_{n+h}= \X_{n+h}-\widehat{\X}_{n+h} $ given $ \X_{n-k+1}, \X_{n-k+2}, \ldots, \X_n$.
\end{enumerate}

Before investigating the theoretical properties of the above bootstrap procedure and evaluating its practical implementation for the construction of prediction bands, some remarks are in order. 

Notice that $ \X^*_{n+h}$ in Step~3 is generated in a model-free way, while the estimated model $ \widehat{g}^*$  is only  used  for obtaining the pseudo-predictor $ \widehat{\X}^*_{n+h}$. In this way the pseudo-error $\X^*_{n+h}-\widehat{\X}_{n+h}^*$ is able to imitate not only the innovation and estimation errors affecting the prediction error  $ \X_{n+h}-\widehat{\X}_{n+h}$ but also the error arising from possible model misspecification.
  In Steps~4 and~5, the backward vector autoregressive representation is used to generate the pseudo time series of scores $ {\bm \xi}_t^*$, $t=1,2, \ldots, n$,  where this pseudo  time series  satisfies the condition $ {\bm \xi}^*_t=\widehat{\bm \xi}_t$ for $ t=n-k+1, n-k+2, \ldots, n$. This enables the generation of a functional pseudo-time series $\X_1^*, \X_2^*, \ldots, \X_n^*$ in Step~6, satisfying requirement~\eqref{eq.requir}.  A problem occurs when $p >k+h$, that is, when the autoregressive order used is larger than the number of future functional observations  needed to 
  run the backward vector autoregression. In this case,  the time series of scores must be extended with  the $ p-(k+h)$ ``missing''  future scores.
  This  problem is solved in Step~3 by generating the additional pseudo-scores $ {\bm \xi}_{n+h+l}^*$ , for $ l=1,2, \ldots, p-(k+h)$.  

\subsection{Bootstrap Validity}\label{sec:3.3}

We establish consistency of the proposed bootstrap procedure in approximating the conditional error distribution~\eqref{eq.prederr} of interest.  Regarding the underlying class of functional processes, we assume that $ {\bf X}$ is a purely non-deterministic, mean square continuous and  $L^4$-${\mathcal M}$ approximable process. The mean square continuity of  ${\bf X}$ implies that its mean and covariance functions are continuous. For simplicity of notation, we assume that $ E\X_t=0$. 

Because we condition on the last $k$ observations, in what follows, all asymptotic results are derived under the assumption that we have observed a functional time series $ \X_s, \X_{s+1}, \ldots, \X_n$ in which we view $n$ as fixed and allow $s\rightarrow -\infty$. This is also the meaning of the statement "as $n\rightarrow \infty$ " used in all derivations and asymptotic considerations in the sequel. Some conditions regarding the underlying process ${\bf X}$ and the behavior of the bootstrap parameters $m$ and $p$ as well as the estimators $\widehat{g}$ and $ \widehat{g}^\ast$ used are first imposed. Notice that to achieve bootstrap consistency, it is necessary to allow for the order $p$ of the fitted autoregression and the dimension $m$ of the number of principal components used to increase to infinity with the sample size. This is required in order for  the bootstrap to appropriately capture both the entire temporal dependence structure of the vector process of scores and the infinite-dimensional structure of the prediction error $ {\mathcal E}_{n+1}$.

{\bf Assumption 1:}
\begin{enumerate}
\item[]
\begin{enumerate}
\item[(i)]   The autocovariance operator $ C_h$ of $ {\bf X}$ satisfies $  \sum_{h\in \mathbb{Z}} |h| \|C_h\|_{HS} <\infty$.
\item[(ii)]      For all $\omega \in [0,\pi]$,  the spectral density operator  ${\mathcal F}_\omega$ is of full rank, that is, $ kern({\mathcal F}_\omega) =0$ and the eigenvalues $\lambda_j$ of the full rank covariance operator $ {\mathcal C}_0$ (in descending order) are   denoted by  $ \lambda_1>\lambda_2> \lambda_3 > \ldots >0 $.
\end{enumerate}
\end{enumerate}

\vspace*{0.2cm}

{\bf Assumption 2:}
The sequences $ p=p(n)$  and $m=m(n)$ satisfy $ p\rightarrow \infty$  and $ m\rightarrow \infty$, as $ n\rightarrow \infty$, such that  
\begin{enumerate}
\item[(i)] \ $ m^2/\sqrt{p} \rightarrow 0$,
\item[(ii)]  \ $ \frac{\displaystyle p^3}{\displaystyle \sqrt{n m} \lambda_m^2}\sqrt{\sum_{j=1}^m\alpha_{j}^{-2}} = O(1)$, where 
$ \alpha_1=\lambda_1-\lambda_2$ and $ \alpha_j=\min\{\lambda_{j-1}-\lambda_j,\lambda_j-\lambda_{j+1}\}$ for $ j=2,3, \ldots, m$.
\item[(iii)] \  $m^4p^2  \|\widetilde{A}_{p,m} - A_{p,m}\| _F =O_P(1)$, where 
$\widetilde{A}_{p,m}=(\widetilde{A}_{1,p},\widetilde{A}_{2,p}, \ldots, \widetilde{A}_{p,p} )$, 
and $  A_{p,m}=(A_{1,p}, $ $  A_{2,p}, \ldots, A_{p,p})$.
Here   $ \widetilde{A}_{j,p}, j=1,2, \ldots, p$ are the same estimators as $ \widehat{A}_{j,p}, j=1,2, \ldots, p$, but based 
on the time series of true scores $ \bm{\xi}_1, \bm{\xi}_2, \ldots, \bm{\xi}_n$. Furthermore,    $ A_{j,p},$ $j=1,2,  \ldots, p$ are the coefficient matrices of the best (in the mean square sense)  linear predictor 
of $ \bm{\xi}_t$ based on the finite past $ \bm{\xi}_{t-1}, \bm{\xi}_{t-2}, \ldots, \bm{\xi}_{t-p}$. 
\end{enumerate} 

{\bf Assumption 3:}  The estimators $\widehat{g} $ and $ \widehat{g}^* $ converge to the same  limit $g_0$; that is, $ \|\widehat{g} -g_0\|_{\mathcal L}=o_P(1)$ and $ \|\widehat{g}^* - g_0\|_{\mathcal L}=o_P(1)$.

Several  comments regarding the above assumptions are in order. Assumption 1(i) implies that the spectral density operator $ {\mathcal F}_\omega$ is a continuously differentiable function of the frequency $\omega$. Regarding Assumption 2, notice first that allowing for the number $m$ of principal components used as well as the order $p$ of the vector autoregression fitted to increase to infinity with the sample size, make the asymptotic analysis of the bootstrap quite involved. This is so because the bootstrap procedure is based on the time series of  estimated instead of  true scores,  the dimension and the order of the fitted vector autoregression increase to infinity, and, at the same time, the eigenvalue $ \lambda_m$ of the lag zero covariance operator $ {\mathcal C}_0$, approaches zero as $ m $ increases to infinity with $n$. As we will see,   a slow increase of $m$ and $p$ with respect to $n$ is required to balance these different effects.

To elaborate, parts  (i)  and  (ii) of Assumption~2  summarize   the  conditions imposed  on  the rate of increase of  $m $ and $p$   to establish 
bootstrap   consistency. Before discussing these conditions in more detail, observe that Assumption~ 2(iii) is a condition that  the estimators of the autoregressive coefficient matrices   have to fulfill, after ignoring the effects  caused by the  fact 
  that  estimated instead of true  scores are used. Observe first  that  in contrast to the estimators $ \widehat{A}_{j,p}$, $j=1,2, \ldots, p$,  based on the vector of estimated scores $ \widehat{\bm \xi}_t$, the estimators 
$ \widetilde{A}_{j,p}$, $j=1,2,\ldots, p$, stated  in Assumption 2(iii) are   based on the true (i.e.,  unobserved) vector of scores  $\bm{\xi}_t$, $t=1,2, \ldots, n$.  As an  an example,  consider    the case where $\widetilde{A}_{j,p}$, $ j=1,2, \ldots, p$, are the well-known  Yule-Walker estimators of $ A_{j,p}$, $j=1,2, \ldots, p$.  By the arguments given in \citet[][p. 3521]{Paparoditis18}, we  have in this case that $ \|\widetilde{A}_{p,m} - A_{p,m}\|_F = O_P\big(mp(\sqrt{m}\lambda_m^{-1} +p)^2/\sqrt{n}\big)$.  From this bound  it is  easily seen by straightforward calculations, that,   for these estimators,  Assumption 2(iii) is satisfied if $ m$ and $ p$ increase to infinity as $ n \rightarrow\infty$ slowly enough such  that $ m^6p^4 = O(\lambda_m^2 \sqrt{n})$ and $ p^2/m^2 =O(\sqrt{n})$.

To give an example of a functional process,  of  an estimator  $ \widehat{A}_j$, $j=1,2, \ldots, p$ and of the    rates  with  which  $ p$ and $m$ have to increase to infinity   so that  all parts of  Assumption 2 are fulfilled, suppose  again that   Yule Walker estimators of $ A_{j,p}$, $j=1,2, \ldots,p$,  are used in the bootstrap procedure. Assume further that  the eigenvalues 
of the lag zero covariance operator  $ {\mathcal C}_0$ satisfy $ \lambda_{j}-\lambda_{j+1} \geq C\cdot j ^{-\vartheta}$, $j=1,2, \ldots$, where $ C$ is a positive constant and $ \vartheta >1$. That is,  assume that  
 the eigenvalues of the lag zero autocovariance operator  $ {\mathcal C}_0$ converge  at a polynomial rate to zero. As it is shown in the supplementary material, Assumption 2(i), (ii) and (iii) are  then satisfied if $ p=O(n^\gamma)$ and $ m = O(n^\delta)$  with $ \gamma>0$ and $\delta>0$,  such  that  
\begin{equation} \label{eq.Ass2Total}
 \gamma \in (0,1/8) \ \mbox{and} \  \delta \in (0,\delta_{max}), \ \mbox{where} \   \delta_{max} = \min\left\{ \frac{1-6\gamma}{6\vartheta}, \ \frac{1-8\gamma}{12+4\vartheta},\  \gamma/4 \right\}.  
\end{equation}
More specifically  and  if  for instance,  $ \vartheta =2$, then   $ (1-8\gamma)/(12+4\vartheta) < (1-6\gamma)/6\vartheta$ and Assumption 2 is satisfied if \[ \gamma \in (0,1/8) \ \ \mbox{and} \ \  \delta \in \big(0, \min\left\{\frac{1-8\gamma}{20}, \gamma/4\right\}\big).\]
Notice that $ 0<\gamma <1/8$ ensures   that $ (1-8\gamma)/20 >0$.

Concerning Assumption 3 and given that we do not focus on a specific predictor, this assumption is necessarily a high-level type assumption. It requires that  the estimator $\widehat{g}^*$, which is based on the bootstrap pseudo-time series $ \X^\ast_1, \X_2^\ast, \ldots, \X_n^\ast$, converges in probability and in operator norm, to the same limit $g_0$ as the estimator $ \widehat{g}$ based on the time series $ \X_1, \X_2, \ldots, \X_n$. Notice that Assumption~3 can only be verified in a case-by-case investigation for a  specific operator $g$ at hand and for the particular estimators $\widehat{g} $ and $\widehat{g}^* $ used to perform the prediction. 
 
To elaborate, consider the following  example.  Suppose  that a  FAR(1) model $ \X_t=\Phi(\X_{t-1}) + \bm{\varepsilon}_t$ is used  in equations~\eqref{eq.xpred} and~\eqref{eq.pred-boot} to obtain the predictors $ \widehat{\X}_{n+1}$ and $ \widehat{\X}^*_{n+1}$,  respectively. A common estimator of $ \Phi $  based on an approximative solution of the Yule-Walker-type equation
 $ C_1=\Phi C_0$, is  given by
 \begin{equation} \label{eq.Phiest}
  \widehat{\Phi}_M(\cdot) = \frac{1}{n-1}\sum_{t=1}^{n-1}\sum_{i=1}^M\sum_{j=1}^M \frac{1}{\widehat{\lambda_j}} \langle\cdot , \widehat{v}_j\rangle
  \langle \X_t , \widehat{v}_j\rangle \langle\X_{t+1}, \widehat{v}_i\rangle \widehat{v}_i,
  \end{equation}
where $M$ is some integer referring to the number of functional principal components included in the estimation of $\Phi$ \citep{Bosq00, HK12}. Observe that~\eqref{eq.Phiest} is a kernel operator with kernel
\begin{equation} \label{eq.kern1}
\widehat{\varphi}_M(\tau,\sigma) = \frac{1}{n-1} \sum_{t=1}^{n-1}\sum_{i=1}^M\sum_{j=1}^M \frac{1}{\widehat{\lambda_j}} 
  \langle \X_t , \widehat{v}_j\rangle \langle\X_{t+1}, \widehat{v}_i\rangle \widehat{v}_{j}(\sigma)\widehat{v}_i(\tau)
\end{equation}
and notice  that $ \widehat{g}=\widehat{\Phi}_M$  in  this example. Furthermore, for  fixed $M$ and by the consistency properties of $ \widehat{\lambda}_j$ and $ \widehat{v}_j$,   it is not difficult to show  that $\| \widehat{\Phi}_M(\cdot)  - g_0\|_{\mathcal L} \stackrel{P}{\rightarrow} 0$,    where  the limiting operator $g_0$ is given by 
\begin{equation} \label{eq.g0M}    
g_0(\cdot) \equiv  C_{1,M} \Big( \sum_{j=1}^M  \frac{1}{\lambda_j}\langle \cdot , v_j\rangle v_j \Big).
\end{equation}
Here, $ C_{1,M}(\cdot) = E \langle \X_{t,M}, \cdot \rangle \X_{t+1,M}$ is a finite rank approximation of the  lag-1 autocovariance operator $ C_1$ (see equation~\eqref{eq.dec} for the definition of $\X_{t,M}$). Further, $ \sum_{j=1}^M  \lambda_j^{-1}\langle \cdot , v_j\rangle v_j $  is  the corresponding approximation of the inverse operator $ C_0^{-1}(\cdot) =  \sum_{j=1}^\infty  \lambda_j^{-1} \langle \cdot , v_j\rangle v_j $, which appears when solving the aforementioned Yule-Walker-type equation \citep[see][Chapter 13, for details]{HK12book}. Similarly,  the same convergence  also holds  true for the bootstrap estimator $\widehat{g}_0^\ast$, that is, $ \|\widehat{g}^\ast_0-g_0\|_{\mathcal L} \stackrel{P}{\rightarrow} 0$ with $ g_0$  given in~\eqref{eq.g0M}. Now,  if interest is focused on consistently estimating the operator $ \Phi$,  then,  from an asymptotic perspective,  the number $M$ of functional principal components used in approximating the inverse of the operator $C_0$, has to increase to infinity at an appropriate rate,  as $n$ goes to infinity.  In this case,  it can be shown 
  under certain regularity conditions that  $ \| \widehat{\Phi}_M - \Phi\|_{\mathcal L} =o_P(1)$ \citep[see][Theorem 8.7]{Bosq00}. Here
$g_0=\Phi$, and this limit is different from  the one  given in~\eqref{eq.g0M}. In such a  case, and for the estimator $\widehat{g}^*$ to also converge to the same limit,  additional arguments  are needed  since the technical derivations are then much more involved,  compared to those used in   the case of a fixed $M$ \citep[we refer to][for more details on this type of asymptotic considerations]{Paparoditis18}. 
 

%
%

Before stating our first consistency  result,  we fix some additional notation. Recall the definition of  $ {\mathcal X}_{n,k}$ and denote by  $ {\mathcal C}_{\mathcal E,h}$ and  $ {\mathcal C}^\ast_{\mathcal E,h}$ the conditional covariance operators  of the random elements $ {\mathcal E}_{n+h}$ and ${\mathcal E}^\ast_{n+h}$, respectively, given $ \X_{n,k}$. That is, $ {\mathcal C}_{\mathcal E,h}=E\big( {\mathcal E}_{n+h}\otimes {\mathcal E}_{n+h}|\X_{n,k}\big)$ and $  {\mathcal C}^\ast_{\mathcal E,h}=E^\ast \big({\mathcal E}^\ast_{n+h}\otimes {\mathcal E}^\ast_{n+h}|\X_{n,k}\big)$, where $ E^\ast$ denotes expectation with respect to the bootstrap distribution. Recall that $ \X_t^*=\X_t$ for $ t\in\{n,n-1, \ldots, n-k+1\}$. Let further,  
\[  \sigma^2_{n+h}(\tau)=c_{\mathcal E,h}(\tau,\tau) \ \ \mbox{and} \ \   \sigma^{\ast^2}_{n+h}(\tau)=c^\ast_{\mathcal E,h}(\tau,\tau), \ \ \tau\in[0,1],\] where $c_{\mathcal E,h} $ and $c_{\mathcal E,h}^\ast $ denote the  kernels of the conditional covariance  (integral) operators $ {\mathcal C}_{\mathcal E,h} $ and $  {\mathcal C}^\ast_{\mathcal E,h}$, respectively, which exist since these operators  are Hilbert-Schmidt.  Denote by $ {\mathcal L}_{\X_{n,k}}(\mathcal E_{n+h})$ the conditional distribution  $ \mathcal E_{n+h}  \big|   \X_{n,k}$, and by $ {\mathcal L}_{\X_{n,k}}(\mathcal E^*_{n+h} |\X_1, \X_2, \ldots, \X_n)$ the conditional distribution $ \mathcal E^*_{n+h} | \X_{n,k} $, given the observed functional time series $ \X_1, \X_2, \ldots, \X_n$. The following theorem establishes consistency of the bootstrap procedure in estimating the conditional distribution of interest.

\begin{theorem} \label{th.boot-val}
Suppose that Assumptions 1, 2, and 3 are satisfied.  Then, 
\begin{equation} \label{eq.weakconvBoot}
d \Big( {\mathcal L}_{\X_{n,k}}(\mathcal E_{n+h}),  {\mathcal L}_{\X_{n,k}}(\mathcal E^*_{n+h} \big| \X_1, \X_2, \ldots, \X_n\big)\Big) = o_P(1),
\end{equation}
where $d$ is any metric metricizing weak convergence on ${\mathcal H}$. 
\end{theorem}

The above result  allows for the use of the conditional distribution of $ {\mathcal E}^*_{n+h} $  to construct pointwise prediction intervals  for $ \X_{n+h}$. Alternatively, consider the use of the conditional  distribution of $ \sup_{\tau\in[0,1]}|{\mathcal E}_{n+h}^*(\tau)|$ to construct prediction bands for $ \X_{n+h}$. Notice that the latter  prediction bands will have the same width for all values of $\tau\in[0,1]$ since they do not appropriately reflect  the local variability of the prediction error $ {\mathcal E}_{n+h}(\tau)$. One way to take  the (possible different) prediction uncertainty  at every $ \tau \in [0,1]$ into account,  is to use the studentized  conditional distribution of the  prediction error, that is to use the process  $\{{\mathcal E}_{n+h}(\tau) / \sigma_{n+h}(\tau), \ \tau \in [0,1]\}$  on ${\mathcal H}$ in order to   construct the prediction bands. However, in this case, and additional to the weak convergence of $ {\mathcal E}_{n+h}^\ast$ to $ {\mathcal E}_{n+h}$ on $ {\mathcal H}$, establishing  bootstrap consistency requires the uniform (over $\tau \in [0,1]$) convergence  of the  conditional variance of the prediction error $ \sigma^{\ast^2}_{n+h}(\tau) $  against  $ \sigma^2_{n+h}(\tau) $. This will allow for the proposed bootstrap procedure to appropriately approximate the random behavior of the studentized process $ \{ {\mathcal E}_{n+h}(\tau)/\sigma_{n+h}(\tau), \tau \in [0,1]\}$. To achieve such a  uniform consistency  of  bootstrap estimates,  additional conditions compared to those stated in the previous Assumptions 2 and 3 are needed.  
We begin with the  following modification of Assumption 2.

\vspace*{0.2cm}

{\bf Assumption 2$^{'}$:} The sequences $ m=m(n)$ and $ p=p(n)$ satisfy  Assumptions 2 (i), (iii), and 
\begin{enumerate}
\item[(ii)] \   $ \frac{\displaystyle p^5\,m}{\displaystyle n^{1/2} \lambda_m^{5/2}}\sqrt{\sum_{j=1}^m\alpha_{j}^{-2}} = O(1)$.
\end{enumerate} 

\vspace*{0.2cm}

Our next assumption imposes additional conditions to those made in Assumption 3 and concern   the mean square consistency properties of the estimators $ \widehat{g}_{(h)}$ and $ \widehat{g}^*_{(h)}$ used to perform the prediction. 

\vspace*{0.2cm}

{\bf Assumption 3$^{'}$:} The estimators $\widehat{g}_{(h)} (x)$ and $ \widehat{g}_{(h)}^*(x) $  satisfy for any given $ x\in {\mathcal H}^k$,\\
 $ \sup_{\tau\in[0,1]}E|  \widehat{g}_{(h)}(x)(\tau)-g_{0,h}(x)(\tau)|^2 \rightarrow 0 $ and $  \sup_{\tau\in[0,1]}E^\ast |  \widehat{g}_{(h)}^*(x)(\tau)-g_{0,h}(x)(\tau)|^2 \rightarrow 0$ in probability.

The following proposition discusses  the conditions that the initial estimators $ \widehat{g}$ and $ \widehat{g}^\ast$ have to fulfill so that  Assumption 3$^{'}$ is satisfied for the important case where the limiting operator $g_0$ is an integral operator. Recall that if $g_0$ is an integral operator with  kernel $c_{g_0}$  satisfying   $ \int_0^1\int_0^1 |c_{g_0}(\tau,s)|^2\tau ds <\infty$,  then $g_0$ is   also  is a Hilbert-Schmidt operator. 

\begin{prop} \label{prop.IO}
Suppose that $ g_0$ is an integral operator with kernel $ c_g :[0,1]\times [0,1] \rightarrow {\mathbb R}$ and let $ \widetilde{c}_g$ be an estimator of  $ c_g$  and $ \widetilde{g}$ the  corresponding integral operator. If $ \|\widetilde{g}\|^2_{HS} \leq  C$ for some constant $ C>0$ and if 
\begin{enumerate}
\item[] (i) \ $E\|\widetilde{g}-g_0\|_{HS}^2 \rightarrow 0$, and, 
\item[] (ii) \ $ \sup_{\tau\in[0,1]} E\int_0^1 \big(\widetilde{c}_g(\tau,s) -c_g(\tau,s)\big)^2ds \rightarrow 0$,
\end{enumerate} 
as $ n\rightarrow \infty$, then Assumption 3$^{'}$ is satisfied for any $h \in {\mathbb N}$. 
\end{prop}

We observe that apart from the basic requirement (i) on the mean square consistency of the estimator $\widetilde{g}$  with respect to the  Hilbert-Schmidt norm, the additional property one needs in case of integral operators is the uniform, mean square consistency  stated in part (ii) of  the above proposition.  \citet[][Theorem 8.7]{Bosq00} established mean square consistency  results when $\widetilde{g}=\widehat{\Phi}_M$ with $ \widehat{\Phi}_M$ given in~\eqref{eq.Phiest},  in case $ g$ is an autoregressive, Hilbert-Schmidt  operator of a FAR(1) process. Recall that  in this case $ \|g\|^2_{HS}= \int_0^1\int_0^1 |c_{\widetilde{g}}(\tau,s)|^2d\tau ds <1$   is required  in order to ensure stationarity and causality of the FAR(1) process. Thus the requirement $ \|\widetilde{g}\|^2_{HS} \leq C$ stated in above  proposition essentially means in this case, that the Hilbert-Schmidt norm of the estimator $ \widetilde{g}$ used, should be bounded away from unity, uniformly in $n$.

We now establish the next theorem, which  concerns the weak convergence of $ {\mathcal L}_{\X_{n,k}}({\mathcal E}_{n+h}^*) $  as well as the uniform convergence of the conditional variance function $ \sigma_{n+h}^{*^2}(\cdot)$ of the bootstrap prediction error. 

\vspace*{0.2cm}

\begin{theorem} \label{th.boot-variance}
Suppose that Assumption 1, 2$^{'}$  and 3$^{'}$ are satisfied. Then, additional to assertion~\eqref{eq.weakconvBoot} of Theorem~\ref{th.boot-val},  the following also holds   true:
\begin{equation} \label{eq.unifvar}
 \sup_{\tau\in[0,1]}\Big| \sigma^{\ast^2}_{n+h}(\tau)-\sigma^2_{n+h}(\tau) \Big| \rightarrow 0, \ \ \mbox{in probability.}
 \end{equation}
\end{theorem}

Theorem~\ref{th.boot-variance} and Slutsky's theorem theoretically justify the use of   $\{{\mathcal E}^*_{n+h}(\tau)$ $/ \sigma^\ast_{n+h}(\tau),  \tau \in [0,1]\}$ to approximate the behavior of  $\{{\mathcal E}_{n+h}(\tau) / \sigma_{n+h}(\tau)) ,  \tau \in [0,1]\}$. As the following corollary shows,  the bootstrap can then  successfully be applied  to construct a  simultaneous prediction band for $ \X_{n+h}$ that appropriately takes into account  the local uncertainty  of prediction. 

\begin{corollary} \label{co.boot-val}
Suppose that the assumptions of Theorem~\ref{th.boot-variance}   are satisfied. For  $\tau\in [0,1]$, let   
\begin{equation*}
V_{n+h}(\tau)=\frac{\X_{n+h}(\tau)-\widehat{\X}_{n+h}(\tau)}{\sigma_{n+h}(\tau)},  \ \  \mbox{and} \ \
V^{*}_{n+h}(\tau)=\frac{\X_{n+h}^*(\tau)-\widehat{\X}^{*}_{n+h}(\tau)}{\sigma_{n+h}^*(\tau)}.\]
Then,
\[ \sup_{x\in\mathbb{R}} \Big|P\Big( \sup_{\tau\in[0,1]} \big| V_{n+h}(\tau) \big| \leq x \big|{\mathcal X}_{n,k} \Big)  - P^\ast\Big( \sup_{\tau\in[0,1]} \big| V^*_{n+h}(\tau) \big| \leq x \big| {\mathcal X}_{n,k}\Big)\Big| \rightarrow 0,\] in probability, where $ P^\ast(A)$ denotes the probability of the event $A$ given the functional  time series $ \X_1, \X_2, \ldots, \X_n$.
\end{corollary}

\section{Practical Construction of Prediction Intervals}\label{sec:4}

As mentioned, the theoretical results of the previous section allow for the use of the quantiles of the distribution of $ \mathcal E^*_{n+h}(\tau)$,  or of  $ V^*_{n+h}(\tau)$, to construct either pointwise prediction intervals for $ \X_{n+h}(\tau)$ for any $\tau \in [0,1]$,  or simultaneous prediction bands for $\{\X_{n+h}(\tau), \tau \in [a,b]\}$ for any  $ 0\leq a < b \leq 1$. Notice that the conditional distributions of ${\mathcal E}^*_{n+h}(\tau)$  and $ V^\ast_{n+h}(\tau)$ can be evaluated by Monte Carlo, that is, by generating $B$ replicates of ${\mathcal E}^*_{n+h}$ and $ \sigma_{n+h}^\ast$,  say, $ \mathcal E^*_{n+h,1}, \mathcal E^*_{n+h,2}, \dots, \mathcal E^*_{n+h,B}$ and $ \sigma_{n+h,1}^\ast, \sigma_{n+h,2}^\ast, \ldots, \sigma_{n+h,B}^\ast$. 
Let $M_{n,h}^* = \sup_{\tau\in [a,b]}\left|V_{n+h}^*(\tau)\right|$, where $V_{n+h}^*(\tau)=\mathcal E^*_{n+h}(\tau)/\sigma^\ast_{n+h}(\tau)$, and denote by $Q_{h,1-\alpha}^*$ the  $1-\alpha$ quantile of the distribution of $M^*_{n,h}$. This distribution   can consistently be estimated using the  $B$ replicates,  $V^\ast_{n+h,b}(\tau)=\mathcal E^*_{n+h,b}(\tau)/\sigma^\ast_{n+h,b}(\tau)$,   $b=1,2, \ldots, B$. The  simultaneous $(1-\alpha)100\%$ prediction band  for $\X_{n+h}$ over the desired interval $[a,b]$, associated with the predictor $ \widehat{\X}_{n+h}$,  is then given by   
\begin{equation} \label{eq.PredBand}
\left\{ \big[\widehat{\X}_{n+h}(\tau)- Q_{h,1-\alpha}^* \cdot\sigma_{n+h}^*(\tau),\,  \widehat{\X}_{n+h}(\tau) + Q_{h,1-\alpha}^*\cdot\sigma_{n+h}^*(\tau)\big],    \quad  \tau \in [a,b] \right\}.
\end{equation}
Clearly,  a pointwise prediction interval for  and any $ \tau \in [0,1]$, or a prediction band for the entire interval $ [0,1]$, can be obtained as special cases of (\ref{eq.PredBand}).
By the theoretical results established in Section~\ref{sec:3}, this    prediction band  achieves (asymptotically) the desired coverage probability  $ 1-\alpha$.

\section{Simulations}\label{sec:5}

\subsection{Choice of Tuning Parameters} 



The theory developed in Section~\ref{sec:3} formulates conditions on the rates at which the bootstrap tuning parameters $m$ and $p$ have to increase to infinity with the sample size $n$ such that bootstrap consistency can be established. An optimal choice of these parameters in practice,  which  also is consistent with the theoretical requirements stated in Section~\ref{sec:3.3}, is left as an open problem for future research. In this section, we will discuss  some relatively simple and practical rules to select $m$ and $p$, which we  found to work well in practice.

We first mention that different approaches have been proposed in the literature on how to choose $m$ and that these approaches also  can  be applied in our setting. We mention here, among others, the pseudo-versions of the Akaike information criterion and Bayesian information criterion considered in \citep{YMW05}; the finite prediction error criterion considered in \citep{ANH15}; the eigenvalue ratio tests \citep{AH13} and the generalized variance ratio criterion introduced in \citep{Paparoditis18}. However, and in order to reduce the computational burden in  our simulation experiments, we apply a simple and commonly used rule to select the number $m$ of functional principal components. This parameter is chosen here using the ratio of the variance explained by the $m$ principal components to the total variance of the random element $\X_t$. 
More specifically, $m$ is selected as 
\begin{equation*}
m_{n,\text{Q}} = \argmin_{m\geq 1} \left\{\frac{\sum^m_{j=1}\widehat{\lambda}_j}{\sum^n_{j=1}\widehat{\lambda}_j} \geq Q\right\},
\end{equation*}
where  $\widehat{\lambda}_s$ denotes the $s$\textsuperscript{th} estimated eigenvalue of the sample lag-0 covariance operator  $ \widehat{\mathcal C}_0$, and $Q$ is a pre-determined value,  with $Q=0.85$ being a common choice \citep[see, e.g.,][p.41]{HK12}.  Once the parameter $m$ has been selected, the order $p$ of the fitted VAR model is chosen using a corrected Akaike information criterion \citep{HT93},  that is, by minimizing 
\begin{equation*}
\text{AICC}(p) = n\log\left|\widehat{\Sigma}_{e,p}\right|+\frac{n(nm+pm^2)}{n-m(p+1)-1},
\end{equation*}
over a range of values of $p$. Here $\widehat{\Sigma}_{e,p} = n^{-1}\sum^n_{t=p+1}\widehat{\bf e}_{t,p}\widehat{\bf e}_{t,p}^{\top}$ and $\widehat{\bf e}_{t,p}$ are the residuals obtained by fitting the VAR($p$) model 
to the $m$-dimensional, vector time series of estimated scores,  $ \widehat{\bm{\xi}}_1, \widehat{{\bm \xi}}_2, \ldots, \widehat{{\bm \xi}}_n$; see also Step~3 of the bootstrap algorithm. 

%

\subsection{Simulation Study}

We utilize Monte Carlo methods to investigate the finite sample performance of the proposed bootstrap procedure. In particular, the final goal of our simulation study is to evaluate the interval forecast accuracy of the bootstrap prediction intervals under both regimes, that is, when the model used for prediction coincides with the model generating the data and when this is not the case. To this end,  and in the first part of our  simulation experiment, we always use a FAR(1) model for prediction. At the same time, we consider a data generating process, which allows for the investigation of the behavior of the proposed bootstrap procedure to construct prediction bands under both aforementioned regimes. The  functional time series $ \X_1, \X_2, \ldots, \X_n$ used stems from the process 
\begin{equation} 
\X_t(\tau) = \int^1_0\psi(\tau,s)\X_{t-1}(s)ds + b\cdot\X_{t-2}(\tau) + B_{t}(\tau) + c \cdot B_{t-1}(\tau), \qquad t=1,2,\dots,n, \label{eq:sim}
\end{equation}
where $\psi(\tau,s) = 0.34\exp^{\frac{1}{2}(\tau^2+s^2)}$,  $ \tau \in [0,1]$, and $B_t(\tau)$ are Brownian motions with zero mean and variance $1/(L-1)$ with $ L=n+1$. Notice that $\|\Psi\|_{\mathcal L}\approx 0.5$, where $ \Psi$ is the integral operator associated with the kernel $ \psi$.  Three parameter combinations are considered:  
\begin{description}
\item[Case I:]  $b=c=0$; \ \ \qquad  \  {\bf Case II:}  $ b=0.4$ and $ c=0$; \ \ \qquad \    {\bf Case III:}  $ b=0.4$ and $ c=0.8$.
\end{description} 
Notice that for $b=c=0$, the data are generated by a FAR$(1)$ model, so in this case, the model used for prediction coincides with the model generating the functional time series. In Case II and for $b=0.4$ and $ c=0$, the data generating process follows a FAR$(2)$ model, which is stationary because $ \|\Psi\|_{\mathcal L} + |b| < 1$.  This case  imitates   a regime of model miss-specification. A situation of an even ``heavier" model misspecification is simulated in Case III,  where for $ b=0.4$ and $c=0.8$, the data generating process is a stationary FARMA(2,1) process with a large moving average component.

Four  sample sizes are considered in the simulation study, $n=100$, $200$, $400$ and $800$. Using the first 80\% of the data as the initial training sample, we compute a one-step-ahead prediction interval. Then, we increase the training sample by one and compute the one-step-ahead prediction interval again. This procedure continues until the training sample reaches the sample size. With 20\% of the data as the testing sample, we compute the interval forecast accuracy of the one-step-ahead prediction based on the FAR(1) model used for prediction. We present results evaluating the performance of the bootstrap method for all three cases described above. 

In the  second part of our simulation experiment, we compare the finite sample performance of the bootstrap for constructing simultaneous prediction bands  with that of the procedure proposed as Algorithm 4 in \citet[][Section 5.2]{ANH15}. In this comparison, the full FARMA(2,1) model is used to generate the data, that is, model~\eqref{eq:sim} with $b=0.4$ and $ c=0.8$. At the same time, the one-step-ahead predictor $ \widehat{\X}_{n+1}$ is obtained using the prediction method proposed by the authors in the aforecited paper. See also case 2) in Section 2.1.

\subsection{Evaluation Criteria of the Interval Forecast Accuracy}

To measure the interval forecast accuracy, we consider the coverage probability difference (CPD) between the nominal coverage probability and empirical coverage probability and the interval score criterion of \cite{GR07}. The pointwise and uniform empirical coverage probabilities are defined as
\begin{align*}
\text{Coverage}_{\text{poinwise}} &= 1- \frac{1}{n_{\text{test}}\times J}\sum^{n_{\text{test}}}_{\eta=1}\sum^J_{j=1}\left[\mathds{1}\big\{\X_{\eta}(\tau_j)>\widehat{\X}_{\eta}^{\text{ub}}(\tau_j)\big\} + \mathds{1}\big\{\X_{\eta}(\tau_j)<\widehat{\X}_{\eta}^{\text{lb}}(\tau_j)\big\}\right], \\
\text{Coverage}_{\text{uniform}}  &= 1-\frac{1}{n_{\text{test}}}\sum^{n_{\text{test}}}_{\eta=1}[\mathds{1}\{\X_{\eta}(\tau) > \widehat{\X}_{\eta}^{\text{ub}}(\tau)\} + \mathds{1}\{\X_{\eta}(\tau)<\widehat{\X}_{\eta}^{\text{lb}}(\tau)\}],
\end{align*} 
where $n_{\text{test}}$ denotes the number of curves in the forecasting period, $J$ denotes the number of discretized data points, $\widehat{\X}_\eta^{ub}$ and $\widehat{\X}_\eta^{lb}$ denote the upper and lower bounds of the corresponding prediction interval, and $\mathds{1}\{\cdot\}$ the indicator function. The pointwise and uniform CPDs are defined as 
\begin{align*}
\text{CPD}_{\text{pointwise}} &= \left|\text{Coverage}_{\text{pointwise}} - \text{Nominal coverage}\right| \\
\text{CPD}_{\text{uniform}} & = \left|\text{Coverage}_{\text{uniform}} - \text{Nominal coverage}\right|.
\end{align*}
Clearly, the smaller the CPD value, the better the performance of the forecasting method. 

The mean interval score criterion introduced by \cite{GR07}, is denoted by $\overline{S}_{\alpha}$ and  combines both CPD and half-width of pointwise prediction interval. It is defined as
\begin{align*}
\overline{S}_{\alpha} = \frac{1}{n_{\text{test}}\times J}\sum^{n_{\text{test}}}_{\eta = 1}\sum^J_{j=1}\Big\{\left[\widehat{\X}_{\eta}^{\text{ub}}(\tau_j) - \widehat{\X}_{\eta}^{\text{lb}}(\tau_j)\right] &+ \frac{2}{\alpha}\left[\X_{\eta}(\tau_j) - \widehat{\X}_{\eta}^{\text{ub}}(\tau_j)\right]\mathds{1}\left[\X_{\eta}(\tau_j) > \widehat{\X}_{\eta}^{\text{ub}}(\tau_j)\right] \\
&+ \frac{2}{\alpha}\left[\widehat{\X}_{\eta}^{\text{lb}}(\tau_j) - \X_{\eta}(\tau_j)\right]\mathds{1}\left[\X_{\eta}(\tau_j) < \widehat{\X}_{\eta}^{\text{lb}}(\tau_j)\right]\Big\}, 
\end{align*}
where 
$\alpha$ denotes the level of coverage, customarily $\alpha=0.2$ corresponding to 80\% prediction interval  and $\alpha=0.05$ corresponding to 95\% prediction interval. The optimal interval score is achieved when $\X_{\eta}(\tau_{j})$ lies between $\widehat{\X}_{\eta}^{\text{lb}}(\tau_{j})$ and $\widehat{\X}_{\eta}^{\text{ub}}(\tau_{j})$, with the distance between the upper bound and lower bounds being minimal.

\subsection{Simulation Results}

As already mentioned, in the first part of our simulations, we   use an  estimated FAR(1) model to  perform the prediction. The corresponding  FAR(1) predictor is obtained as  $ \widehat{\X}_{n+1} = \overline{\X}_n + \widehat{\Phi}(\X_{n}- \overline{\X}_n)$, where  $ \widehat{\Phi} $ is  a regularized, Yule-Walker-type estimator; see also~\eqref{eq.Phiest}. Table~\ref{tab1:sim} and Table~\ref{tab2:sim}  present results  based on $1,000$ replications (i.e., a pseudo-random seed for each replication) and $ B=1,000$ bootstrap repetitions.
 For the case  $n=800$ and for computational reasons, we only consider $300$ replications. 
 Table~\ref{tab1:sim} presents  results for Case I  ($b=c=0$) and Case II ($b=0.4, c=0$) while  Table~\ref{tab2:sim} for Case III ($b=0.4,c=0.8$).
\begin{table}[!htbp]
\centering
\tabcolsep 0.07in
\renewcommand{\arraystretch}{0.9} 
\caption{Empirical performance of the bootstrap prediction intervals and bands using the FAR(1) model to perform one-step-ahead predictions for functional time series stemming from model~\eqref{eq:sim} with $c=0$ and for Case I ($b=0$) and Case II ($b=0.4$).} 
\label{tab1:sim}
\begin{tabular}{@{}llrrrrrrrr@{}}
\toprule
Nominal & & \multicolumn{2}{c}{$n=100$} & \multicolumn{2}{c}{$n=200$} & \multicolumn{2}{c}{$n=400$} & \multicolumn{2}{c}{$n=800$} \\
coverage & Criterion & $b=0$ & $b=0.4$ & $b=0$ & $b=0.4$ & $b=0$ & $b=0.4$ & $b=0$ & $b=0.4$ \\ 
\midrule 
80\% & Coverage$_{\text{pointwise}}$ & 0.778 & 0.745 & 0.791 & 0.778 & 0.797 & 0.797 & 0.799 & 0.797 \\ 
& CPD$_{\text{pointwise}}$ & 0.0497 & 0.0739 & 0.0344 & 0.0413 & 0.0230 &  0.0238 & 0.0167 & 0.0180 \\
& Coverage$_{\text{uniform}}$ & 0.772 & 0.714 & 0.798 & 0.766 & 0.812 & 0.792 & 0.803 & 0.796 \\
& CPD$_{\text{uniform}}$  & 0.0841 & 0.1152 & 0.0551 & 0.0647 & 0.0390 & 0.0391 & 0.0271 & 0.0283 \\
& $\overline{S}_{\alpha=0.2}$ & 2.5322 & 2.9329 & 2.5346 & 2.7687 & 2.4764 & 2.6647 & 2.4398 &  2.6141 \\
\\
95\% &  Coverage$_{\text{pointwise}}$ & 0.925 & 0.897 & 0.932 & 0.919 & 0.936 & 0.930 & 0.943 & 0.940 \\
& CPD$_{\text{pointwise}}$ & 0.0344 & 0.0580 & 0.0245 & 0.0344 & 0.0172 &  0.0220 & 0.0164 & 0.0202 \\
& Coverage$_{\text{uniform}}$ & 0.920 & 0.877 & 0.932 & 0.909 & 0.941 & 0.927 & 0.946 & 0.934 \\
& CPD$_{\text{uniform}}$ &      0.0516 & 0.0862 & 0.0338 & 0.0502 & 0.0230 &  0.0303 & 0.0185 & 0.0282 \\
& $\overline{S}_{\alpha=0.05}$  &  3.4949 & 4.3379 & 3.5005 & 3.9347 & 3.3988 & 3.7328  & 3.3558 & 3.6646 \\
\bottomrule
\end{tabular}
\end{table}

\begin{table}[!htbp]
\centering
\tabcolsep 0.3in
\renewcommand{\arraystretch}{0.9} 
\caption{Empirical performance of the bootstrap prediction intervals and bands using the FAR(1) model to perform one-step-ahead predictions for functional time series stemming from model~\eqref{eq:sim} with $b=0.4$ and $c=0.8$.} \label{tab2:sim}
\begin{tabular}{@{}llrrrr@{}}
\toprule
 Nominal & & \\
 coverage & Criterion & $n=100$ & $n=200$ & $n=400$  & $n=800$ \\
\midrule
 80\% & Coverage$_{\text{pointwise}}$ & 0.741 & 0.782 & 0.803 & 0.820 \\
 & CPD$_{\text{pointwise}}$ & 0.0888 & 0.0467 & 0.0275 & 0.0235 \\
 & Coverage$_{\text{uniform}}$ & 0.735 & 0.793 & 0.827 & 0.853 \\
 & CPD$_{\text{uniform}}$ & 0.1164 &  0.0652 & 0.0469 & 0.0542 \\
 & $\overline{S}_{\alpha = 0.2}$ & 3.5907 & 3.0588 & 2.8808 & 2.7326 \\
\\
 95\% & Coverage$_{\text{pointwise}}$ & 0.890 & 0.919 & 0.931 & 0.942 \\
 & CPD$_{\text{pointwise}}$ & 0.0702 & 0.0372 & 0.0220 & 0.0114 \\
 & Coverage$_{\text{uniform}}$ & 0.875 & 0.914 & 0.934 & 0.950 \\
 & CPD$_{\text{uniform}}$ &  0.0921 & 0.0484 & 0.0271 & 0.0142 \\
 & $\overline{S}_{\alpha = 0.05}$ & 5.8624 & 4.4679 & 4.0899 & 3.8067 \\
\bottomrule
\end{tabular}
\end{table}

From  Table~\ref{tab1:sim} and Table~\ref{tab2:sim}  some interesting observations can be made. First of all, the empirical coverage of the prediction intervals is good, even for the small sample sizes considered. It improves fast and considerably as the sample size increases, and they get quite close to the desired nominal coverages. This is true for both the pointwise prediction intervals and the simultaneous prediction bands considered and for both coverage levels  used in the simulation study. Further, the $\overline{S}_{\alpha} $ values are systematically larger for the cases   $b=0.4,c=0$ and $b=0.4,c=0.8$ than for the case $b=c=0$. As discussed in the introduction, this expected result is attributable to the fact that the model misspecification errors occurring for  $b=0.4$ and $c=0.8$, respectively, for $b=0.4$ and $c=0$, also  cause an increase in the variability of the prediction error distribution, leading to prediction intervals that are wider than those for $b=0$ and $c=0$. Consequently,  for the ``heavier'' case of model misspecification, Case III, that is for $ b=0.4$ and $ c=0.8$, the $\overline{S}_\alpha$ values are larger than for the case $ b=0.4$ and $c=0$. 

We next present in Table~\ref{table3} results for the second part of our simulations, which compare  the performance of the bootstrap method proposed in this paper with  Algorithm 4 in \cite{ANH15} and used  to construct one-step-ahead prediction bands. Notice that  the predictor $ \widehat{\X}_{n+1}$ used in this part of the simulation study is the one  proposed in \cite{ANH15}. As it is seen from Table~\ref{table3}, the bootstrap approach proposed outperforms the aforementioned algorithm  in \cite{ANH15}, in that the coverage rates are uniformly closer to the nominal levels and the mean interval scores $ \overline{S}_\alpha$, are smaller for all $\alpha$ and for both sample sizes considered.

\begin{table}[!hbp]
\centering
\tabcolsep 0.09in
\caption{Finite sample performance of the bootstrap prediction intervals and of the prediction intervals of \citeauthor{ANH15} \citeyearpar{ANH15} for time series stemming from model~\eqref{eq:sim} with $b=0.4$ and $c=0.8$.}\label{tab:compar_Aue}
\begin{tabular}{@{}llrrrrr@{}}
\toprule
Nominal & & \multicolumn{2}{c}{$n=100$} & \multicolumn{2}{c}{$n=200$} \\
coverage & Criterion & \citeauthor{ANH15}'s \citeyearpar{ANH15} & Bootstrap& \citeauthor{ANH15}'s \citeyearpar{ANH15} & Bootstrap \\
\midrule
80\% & Coverage$_{\text{pointwise}}$ & 0.911 & 0.837 & 0.880 & 0.833 \\
& CPD$_{\text{pointwise}}$ & 0.1140 & 0.0610 & 0.0901 & 0.0348 \\
& Coverage$_{\text{uniform}}$ & 0.915 & 0.842 & 0.884 & 0.824 \\
& CPD$_{\text{uniform}}$ & 0.1255 & 0.0897 & 0.0921 & 0.0466 \\
& $\overline{S}_{\alpha=0.2}$ & 6.0432 & 3.0487 & 5.4047 & 2.7993 \\
\\
95\% & Coverage$_{\text{pointwise}}$ & 0.977 & 0.953 & 0.967 & 0.949 \\
& CPD$_{\text{pointwise}}$ & 0.0340 & 0.0291 & 0.0235 & 0.0106 \\
& Coverage$_{\text{uniform}}$ & 0.974 & 0.946 & 0.968 & 0.949  \\
& CPD$_{\text{uniform}}$ & 0.0400 & 0.0435 & 0.0240 & 0.0154 \\
& $\overline{S}_{\alpha=0.05}$ & 7.1303 & 4.3025 & 5.7686 & 3.8679 \\
\bottomrule \label{table3}
\end{tabular}
\end{table}

\section{Empirical Data Analysis}\label{sec:6}

For the two real-life data sets analyzed in this section and in the supplementary material, we consider  (one and two step head)  prediction using   a FAR(1) model and a nonparametric forecasting  method (NFR). The latter method applied for one-step-ahead prediction, uses a nonparametric estimation  of the lag-1 conditional mean function $g(\X_n)=E(\X_{n+1} | \X_n)$ (also see Section~\ref{sec:2.1}).  Recall that $ g(\X_n)$ is the best (in the mean square sense) predictor of $ \X_{n+1}$  based on $ \X_n$, and that different  data-driven smoothing techniques 
exist for estimation of  $g$. We refer to the functional Nadaraya-Watson estimator \citep[see, e.g.,][]{Masry05, FV06},  the functional local linear estimator, \cite{BEM11}, the functional $k$-nearest neighbor estimator,   \cite{KV13}, and the distance-based local linear estimator \cite{BDF10}, to name a few. We use for the one-step-ahead prediction the Nadaraya-Watson estimator of $g$, which leads to the predictor $\widehat{\X}_{n+1}=\widehat{g}(\X_n)$ with $ \widehat{g}$  given by 
\begin{equation*}
\widehat{g}(x) = \frac{\sum^n_{t = 2}K\left(d(x, \X_{t-1})/\delta\right)\X_{t}}{\sum^n_{t= 2}K\left(d(x, \X_{t-1}))/\delta\right)}. 
\end{equation*}
In the Nadaraya-Watson estimator, $K(\cdot)$ is the Gaussian kernel and $\delta$ is bandwidth, which in our calculations has been obtained using a generalized cross-validation procedure. 
 
In addition to the construction of prediction bands, we also demonstrate how the proposed bootstrap method can be used to select the prediction method that performs better according to  some user-specified criterion. In particular, since  the future random element $\X_{n+1}^\ast$ is generated in a model-free way,  the bootstrap prediction error $ \widehat{\X}_{n+1}^\ast - \X^\ast_{n+1}$ correctly imitates the  behavior of the prediction error $ \widehat{\X}_{n+1}-\X_{n+1}$ associated with  the particular model/method used to perform the prediction. Thus, using some loss $L(\widehat{\X}_{n+1}, \X_{n+1})$ and based on the behavior of the corresponding bootstrap loss $L(\widehat{\X}^\ast_{n+1},\X^\ast_{n+1})$, the proposed bootstrap procedure can also be applied to select between different predictors,  the one 
which  performs better. In the following, we demonstrate such an application of the bootstrap in selecting between the FAR(1) and the NFR method   using  the behavior of the bootstrap estimator of the  (conditional) mean square error of prediction, i.e., of  $ E^\ast[(\widehat{\X}^\ast_{n+1}(\tau)-\X^\ast_{n+1}(\tau))^2|{\mathcal \X}_n ]$, $\tau \in [0,1]$.
 
\subsection{Monthly Sea Surface Temperature Data Set}\label{sec:6.1}

\begin{figure}[!htbp]
\centering
\subfloat[Univariate time series plot]
{\includegraphics[width=8.3cm]{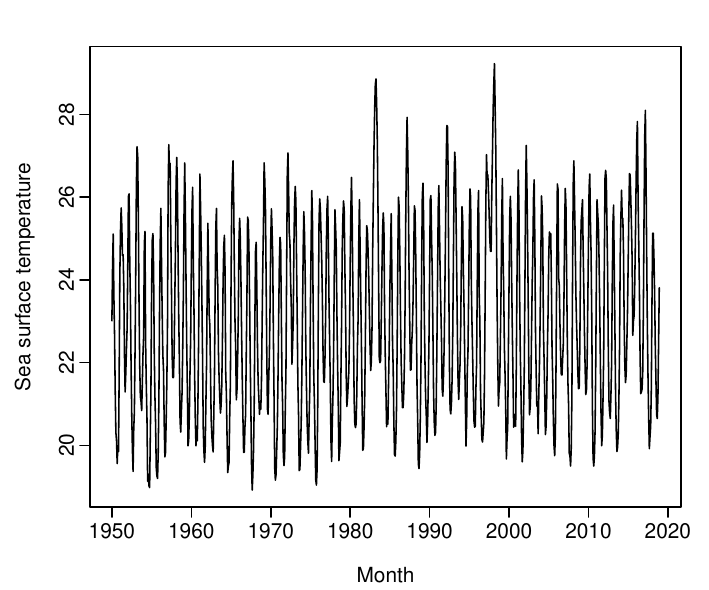}\label{fig:2a}}
\qquad
\subfloat[Functional time series plot]
{\includegraphics[width=8.3cm]{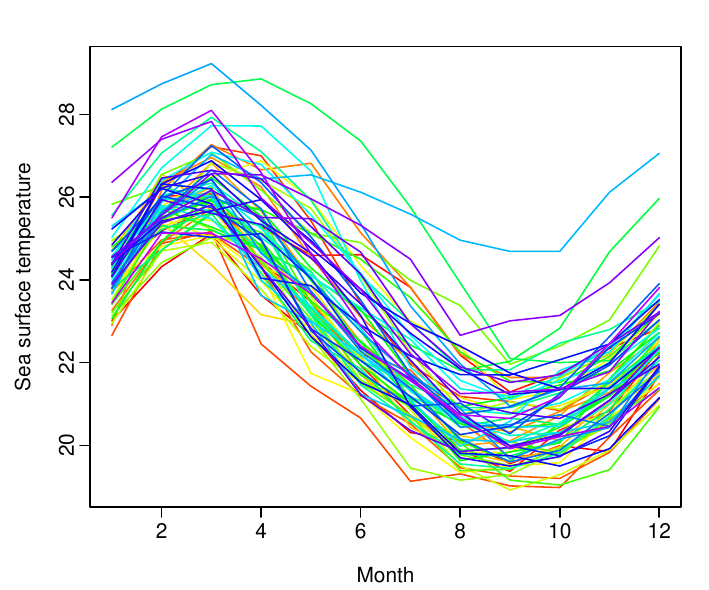}\label{fig:2b}}
\caption{Time series (left panel) and rainbow plots (right panel) of sea surface temperatures in Ni\~{n}o 1+2 region from January 1982 to December 2017.}\label{fig:2}
\end{figure}

Consider the monthly sea surface temperatures from January 1950 to December 2018 available at \url{https://www.cpc.ncep.noaa.gov/data/indices/ersst5.nino.mth.81-10.ascii}. These averaged sea surface temperatures were measured by moored buoys in the ``Ni\~{n}o region". We consider all four Ni\~{n}o regions: Ni\~{n}o 1+2 is defined by the coordinates $0-10^{\circ}$ South, $90-80^{\circ}$ West; Ni\~{n}o 3 is defined by the coordinates $5^{\circ}$ North -- $5^{\circ}$ South, $150^{\circ}$ -- $90^{\circ}$ West; Ni\~{n}o 4 is defined by the coordinates $5^{\circ}$ North -- $5^{\circ}$ South, $160^{\circ}$ East -- $150^{\circ}$ West; Ni\~{n}o 3+4 is defined by the coordinates $5^{\circ}$ North -- $5^{\circ}$ South, $170-120^{\circ}$ West. For the sea surface temperatures in Ni\~{n}o 1+2 region, univariate and functional time series plots are shown in Figure~\ref{fig:2}.

\begin{figure}[!htbp]
\centering
{\includegraphics[width=8.3cm]{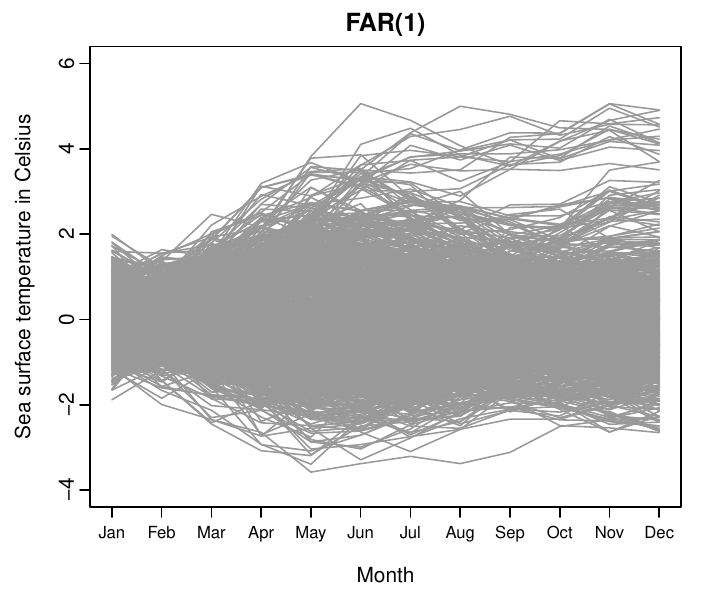}\label{fig:prediction_err}}
\quad
{\includegraphics[width=8.3cm]{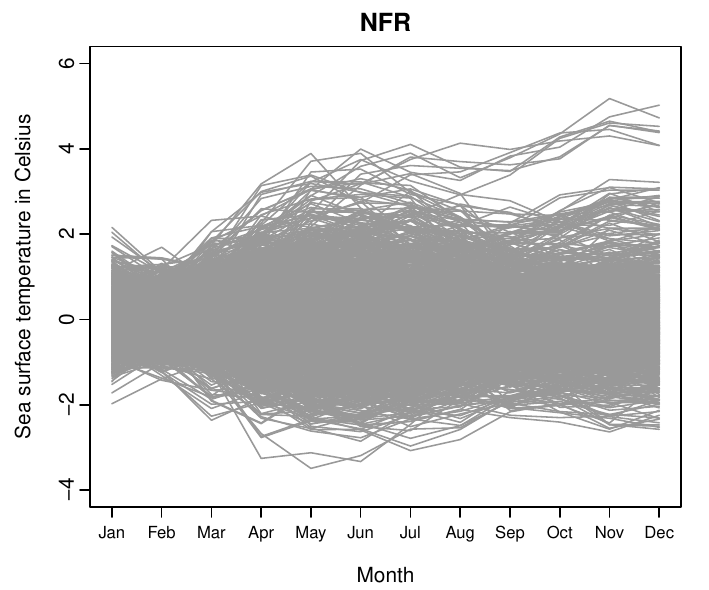}\label{fig:prediction_err_2017}}
\caption{Bootstrap-generated one-step-ahead prediction error curves for the sea surface temperature data using the FAR(1) model  (left panel) and the NFR  method (right panel).}\label{fig:5}
\end{figure}

Applying the proposed bootstrap procedure and the two compared forecasting methods, we generate a set of $B=1,000$ bootstrap one-step-ahead prediction error curves $ \X^\ast_{n+1}-\widehat{\X}^\ast_{n+1}$. These are presented in the left panel of Figure~\ref{fig:5} for the FAR(1) predictor and in the right panel of the same figure for the NFR predictor.
 
\begin{figure}[!htbp]
\centering
\subfloat[Mean square error of one-step-ahead prediction errors]
{\includegraphics[width=8.3cm]{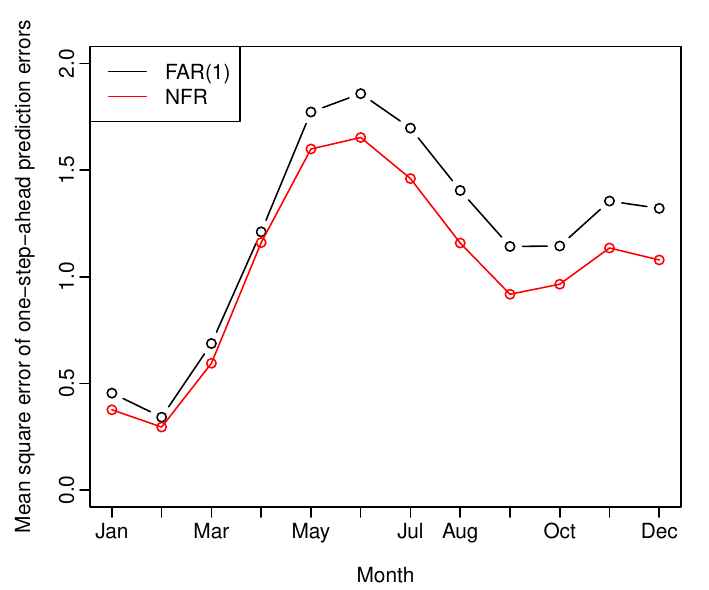}\label{fig:55_a}}
\quad
\subfloat[Pointwise and uniform prediction intervals (NFR method)]
{\includegraphics[width=8.3cm]{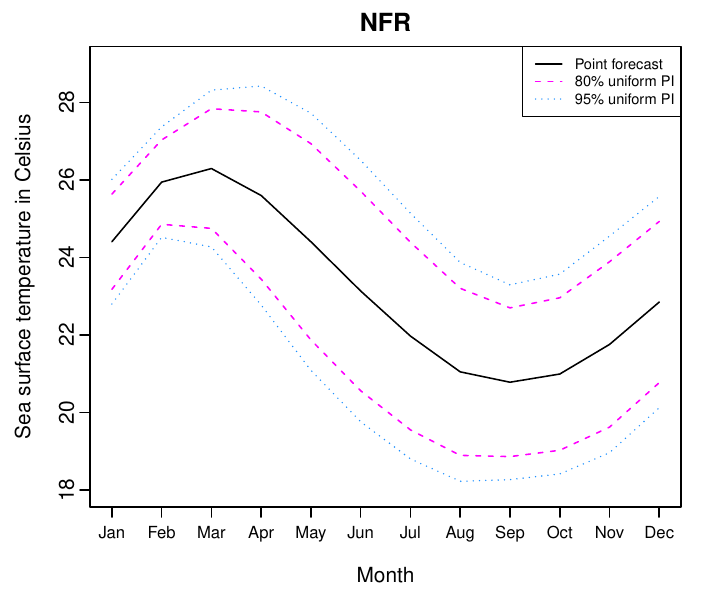}\label{fig:55_b}}
\caption{Bootstrap estimates of the mean square error of the one-step-ahead prediction  for the FAR(1) model and the NFR method (left panel). Point forecast together with 80\% and 95\%   simultaneous  prediction bands using the NFR method (right panel).}\label{fig:55}
\end{figure}

In Figure~\ref{fig:55_a}, we present the bootstrap estimates of the (conditional) mean square error 
$ E[(\widehat{\X}_{n+1}(\tau_j) - \X_{n+1}(\tau_j))^2|\X_n]$, $ j=1,2, \ldots, J$,  of the two prediction methods. As shown in this figure, the mean square error produced by the NFR method is uniformly (across all points $\tau_1, \tau_2, \ldots, \tau_J$) smaller than the corresponding mean square prediction error produced using the FAR(1) method. Thus, for this functional time series, using the NFR method to perform the one-step-ahead prediction seems preferable. The point prediction using this method as well as  the corresponding $ 80\%$ and $ 95\%$ simultaneous prediction bands are shown in Figure~\ref{fig:55_b}.

\section{Conclusions}\label{sec:7}

We have presented a novel bootstrap method for the construction of prediction bands for a functional time series. In a model-free way, our method generates future functional pseudo-random elements that allow for valid estimation of the conditional distribution of the prediction error that a user-selected prediction method produces. The obtained bootstrap estimates of the prediction error distribution consider the innovation and the estimation errors associated with prediction and the error arising from a different model being used for prediction than the one generating the functional time series observed. Theoretical results were presented to justify using the proposed bootstrap method in constructing prediction bands that also appropriately consider the local variability of the conditional distribution of the prediction error. We have demonstrated the good finite sample performance of the bootstrap method presented through a series of simulations. The real-life data sets analyzed have demonstrated the capabilities and the good finite sample behavior of the presented bootstrap method for the construction of prediction bands.

%
%

\bigskip
\begin{center}
{\large\bf SUPPLEMENTARY MATERIAL}
\end{center}

This supplementary material contains the analysis of a second real data example, the derivation of Condition (16), some auxiliary lemmas and 
the proofs of the theoretical results presented in the main paper. 

\section{Intraday PM$_{10}$ Data Set}

\begin{figure}[!htbp]
\centering
\subfloat[A univariate time series display]
{\includegraphics[width=7.7cm]{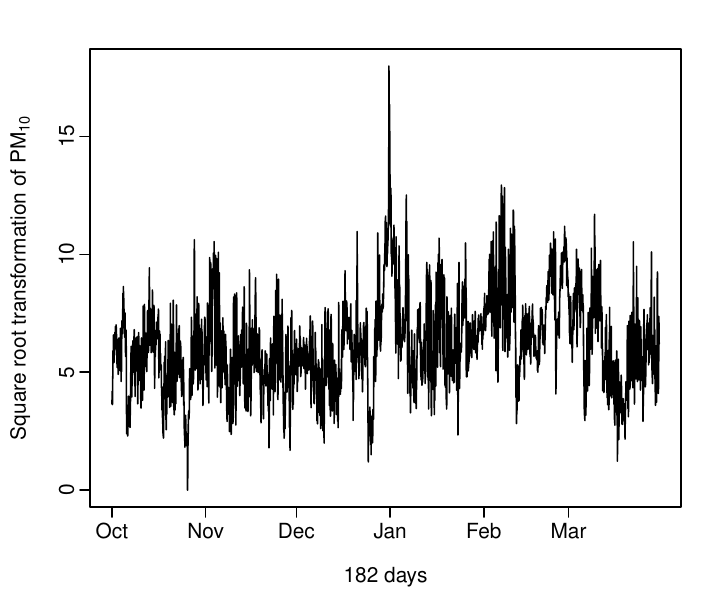}\label{fig:pm_univariate}} 
\qquad
\subfloat[A functional time series display]
{\includegraphics[width=7.7cm]{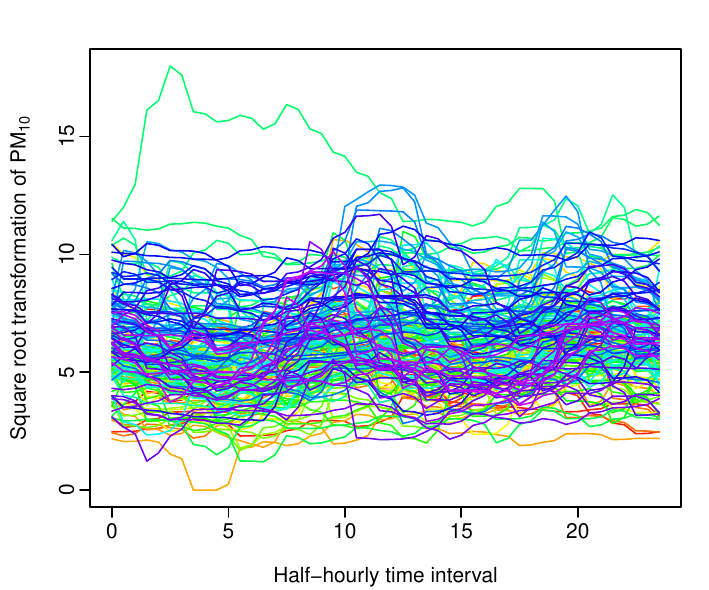}\label{fig:pm_fts}}
\caption{Graphical displays of intraday measurements of the PM$_{10}$ from October 1, 2010, to March 31, 2011, in Graz, Austria.}\label{fig:pm}
\end{figure}

We analyze the half-hourly measurements of the concentration of particulate matter with an aerodynamic diameter of less than 10um in ambient air taken in Graz, Austria, from October 1, 2010, to March 31, 2011. We convert $N=8,736$ discrete univariate time series points into $n=182$ daily curves. A univariate time series display of intraday pollution curves is given in Figure~\ref{fig:pm_univariate}, with the same data shown in Figure~\ref{fig:pm_fts} as a time series of functions.

\begin{figure}[!htbp]
\centering
{\includegraphics[width = 7.6cm]{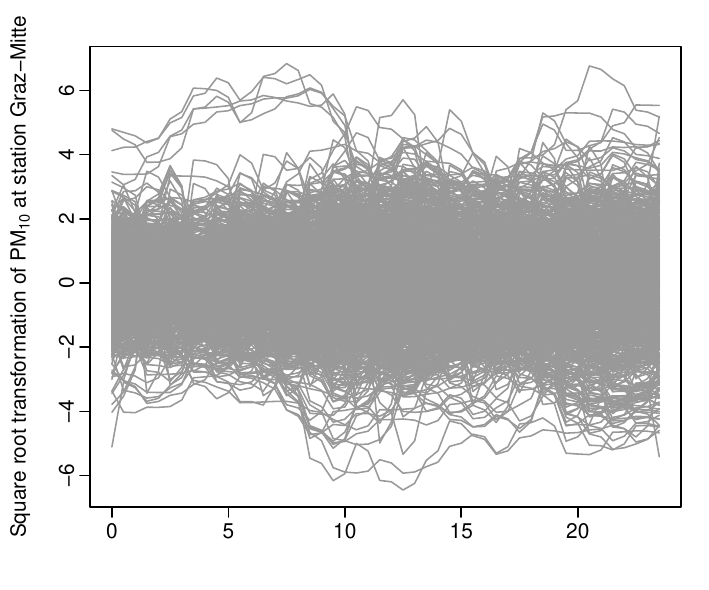}\label{}}
\qquad
{\includegraphics[width = 7.6cm]{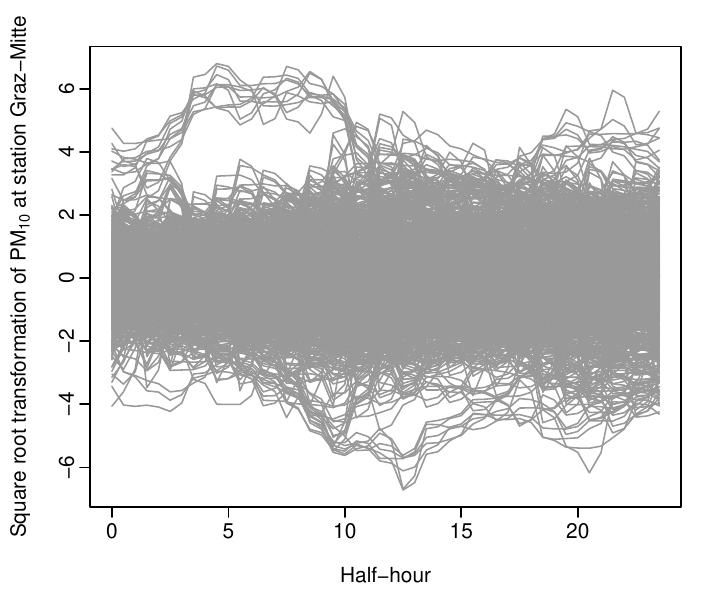}\label{}}
\caption{One-step-ahead prediction error curves for the intraday PM$_{10}$  data using the FAR(1) model (left panel) and the NFR method (right panel).} \label{fig:pm_fore_predict}
\end{figure}

\begin{figure}[!t]
\centering
{\includegraphics[width=7.7cm]{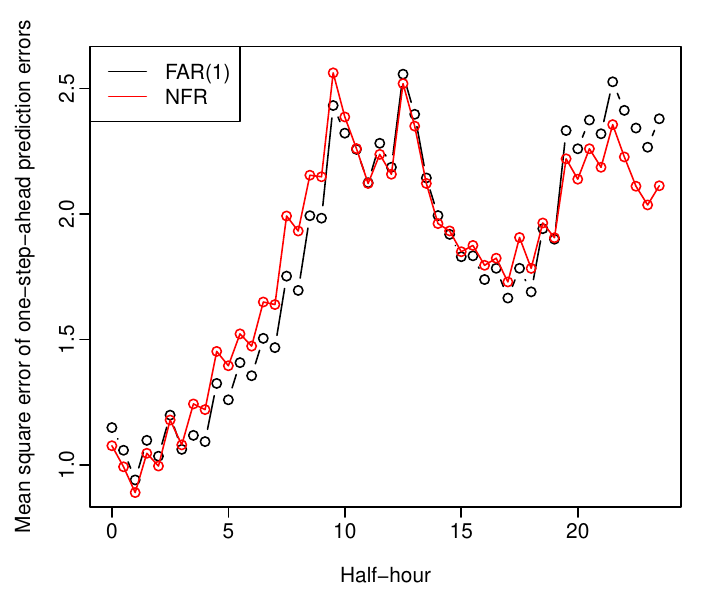}}
\quad
{\includegraphics[width=7.7cm]{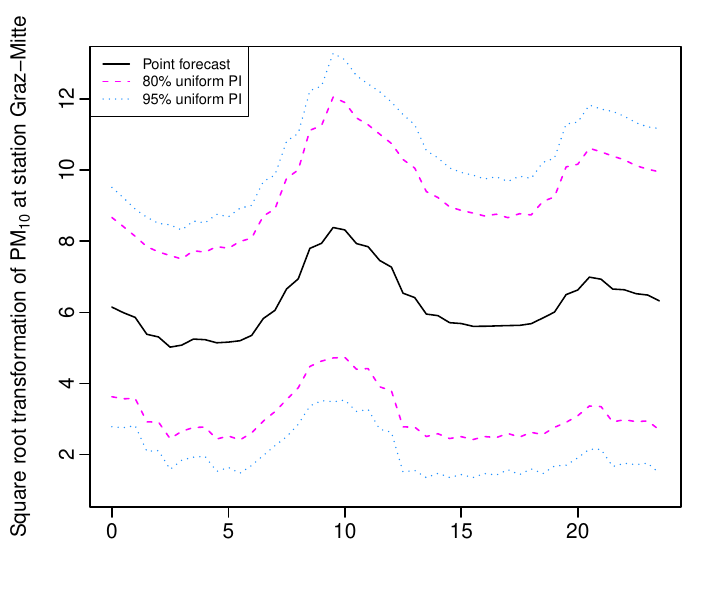}}
\caption{Bootstrap estimates of the mean square error of the one-step-ahead prediction  for the FAR(1) model and the NFR method (left panel). Point forecast together with 80\% and 95\%   simultaneous   prediction bands using the FAR(1) method (right panel).}\label{fig:85}
\end{figure}

Using the bootstrap procedure, we construct one- and two-step-ahead prediction bands for this functional time series. We first generate $B=1,000$ functional pseudo-time series, and we apply the FAR(1) and the NFR forecasting methods to obtain a set of one-step-ahead prediction error curves. These are displayed in Figure~\ref{fig:pm_fore_predict}. In the left panel of Figure~\ref{fig:85}, we show the bootstrap estimates of the (conditional) 
mean square prediction error $ E[(\widehat{\X}_{n+1}(\tau_j) - \X_{n+1}(\tau_j))^2 |\X_n]$, $ j=1,2, \ldots, J$, obtained using the FAR(1) and NFR methods. As this figure shows, neither of the two methods is uniformly (i.e., across all points $\tau_j\in[0,1]$) better, with the FAR(1) method having  a slight advantage. In particular, the bootstrap estimated conditional, root mean square error (RMSE) of the FAR(1) method is 1.451 compared with 2.054 of the NFR method. 

We apply the FAR(1) method to perform the prediction for $h=1$ and $h=2$. The corresponding bootstrap-based simultaneous prediction bands for $h=1$, are displayed in the right panel of Figure~\ref{fig:85}.

In contrast to the relatively small sample size of $n=69$ curves of the Monthly Sea Surface Temperature data analyzed in the previous example, the sample size of $n=182$ curves of the intraday PM$_{10}$ data considered in this section allows us to further evaluate  the performance of the FAR(1) prediction method. For this, we use the observed data from October 1, 2010, to January 30, 2011, as the initial training sample, and we produce one- and two-step-ahead forecasts by increasing the training sample by one curve each time. We iterate this procedure until the training sample contains all the observed data. In this way, we construct 60 one-step and 59 two-step-ahead forecasts, enabling us to assess the interval forecast accuracy of the FAR(1) method.  The results obtained for $h=1$ and $h=2$ are shown in Table~\ref{tab:pm_results}. The FAR(1) prediction performs well for this functional time series, and the proposed bootstrap method produces for both forecasting steps, prediction intervals, the empirical coverages of which are close to the desired nominal levels. Notice that, as it is expected, the uncertainty associated with  $h=2$ is larger than that for  $h=1$. 

\begin{table}[!htbp]
\tabcolsep 0.55in
\centering
\caption{Evaluation of the interval forecast accuracy at the forecast horizons $h=1$ and $h=2$ for the PM$_{10}$ data set using the FAR(1) model for prediction and the sieve bootstrap procedure with 1000 bootstrap replications.}\label{tab:pm_results}
\begin{tabular}{@{}lllrr@{}}
\toprule
Nominal  \\
 coverage & Criterion &  &  $h=1$ &  $h=2$ \\
\midrule
80\% & CPD$_{\text{pointwise}}$&  &   0.789 & 	0.777	\\
& $\overline{S}_{\alpha=0.2}$& & 5.189 &   7.572 \\
	& CPD$_{\text{uniform}}$	&  & 0.833 & 0.777 \\
\\
95\% & CPD$_{\text{pointwise}}$&  & 0.943 &  0.923 \\
	& $\overline{S}_{\alpha=0.05}$	&  & 7.528 &  10.941 \\
	& CPD$_{\text{uniform}}$ 	&  &  0.917 &  0.901 \\
\bottomrule
\end{tabular}	
\end{table}

\section{Auxiliary Results and Proofs}\label{sec:proofs}

%
%
%
%

Recall that $ \X_{n+h}= \sum_{j=1}^m {\bf 1}^\top_{j}\bm{\xi}_{n+h}v_j + U_{n+h,m}$, $\widehat{\X}_{n+h} =  \widehat{g}_{(h)}(\X_n,\ldots, \X_{n-k+1})$, 
$ \X^*_{n+h}=  $   \\  $ \sum_{j=1}^m {\bf 1}^\top_{j}\bm{\xi}^*_{n+h}\widehat{v}_j + U^*_{n+h,m}$ and $\widehat{\X}^*_{n+h} =\widehat{g}^*_{(h)}(\X_n,\ldots, \X_{n-k+1})$. Define $ \X^+_{n+h,m} = \sum_{j=1}^m{\bf 1}_j^\top {\bm \xi}_{n+h}^+ v_j $, where $ {\bm \xi}_{n+h}^+=\sum_{j=1}^p \widetilde{A}_{j,p} {\bm \xi}_{n+h-j} $ $ + {\bm e}^+_{n+h}$, with  ${\bm e}^+_{n+h} $ i.i.d. resampled from the set $\{\widetilde{\bm e}_t -\overline{\widetilde{\bm e}}_n, t=p+1, p+2, \ldots, n \}$, $ \overline{\widetilde{\bm e}}_n=(n-p)^{-1}\sum_{t=p+1}^n \widetilde{\bm e}_t$,
and $\widetilde{\bm e}_t={\bm \xi}_t-\sum_{j=1}^p \widetilde{A}_{j,p} {\bm \xi}_{t-j} $, $t=p+1, p+2, \ldots, n$,  are the residuals obtained from an  autoregressive fit based on the time series of true scores ${\bm \xi}_1, {\bm \xi}_2, \ldots, {\bm \xi}_n$.   

\vspace*{0.3cm}

\noindent {\bf Derivation  of Condition in (16):}  Note that the assumption $ \lambda_j-\lambda_{j+1} \geq Cj^{-\vartheta}$ for $j=1,2, \ldots$,  implies that $1/\lambda_j \leq C^{-1} j^\vartheta$.  For $ p=O(n^\gamma)$ and  $ m=O(n^\delta)$, with $\gamma >0$ and $ \delta >0$,  Assumption 2(i) is satisfied if $ \delta <\gamma/4$. Regarding Assumption 2(ii)  verify that 
\[ \sqrt{\sum_{j=1}^m \frac{1}{\alpha_j^2}  }\leq \frac{1}{C} \frac{1}{\sqrt{2\vartheta +1}} (m+1)^{\vartheta +1/2} = O(n^{\delta (\vartheta +1/2)}),\]
From this  we get that 
\[ \frac{p^3}{\sqrt{nm} \lambda_m^2} \sqrt{\sum_{j=1}^m \frac{1}{\alpha_j^2}}   = O \big( n^{3\gamma -1/2 + 3\vartheta \delta}\big).\]
i.e.,  Assumption 2(iii) is satisfied if $\delta \leq (1-6\gamma)/(6\vartheta)$ and $ \gamma <1/6$.
For Assumption 2(iii)  and in the case of Yule Walker estimators, we have     
\[ \frac{m^6p^4}{\lambda_m^2 \sqrt{n}} = O(n^{6\delta +4\gamma +2\vartheta \delta -1/2}),\]
which is fulfilled for    $ \delta \leq (1 -8\gamma)/(12+4\vartheta)$ and $ \gamma <1/8$ . Note  that 
   for $ \gamma <1/8$ and $ \delta >0$, the condition $p^2/m^2 =O(\sqrt{n})$ also  holds true. Putting   the derived requirements  on $\gamma$ and $\delta$   together, leads to   condition in (15). \hfill $\Box$

\begin{lem} \label{le.papa2018} Let $\Gamma_m(0)=E({\bm \xi}_{n+h}{\bm \xi}_{n+h}^\top)$, 
$\Gamma_m^+(0)=E({\bm \xi}^+_{n+h}{\bm \xi}_{n+h}^{+^\top}) $ and $ \Gamma_m^*(0)=E({\bm \xi}^*_{n+h}{\bm \xi}_{n+h}^{*^\top} )$.  If Assumption 1 and 2 are satisfied, then, 
\[\|\Gamma_m^+(0)-\Gamma_m(0)\|_F = O_P\Big(\frac{m^2}{\sqrt{p}}\Big).\]
If Assumption 1 and 2' are satisfied, then 
\[    \|\Gamma_m^*(0)-\Gamma_m(0)\|_F = O_P\Big(\frac{p^5\sqrt{m}}{\sqrt{n}\lambda^2_m}\sqrt{\sum_{j=1}^m\alpha_j^{-2}} \Big).\]
\end{lem}

{\bf Proof:} Let $\Psi_{j,p}$, $\widetilde{\Psi}_{j,p}$ and $\widehat{\Psi}_{j,p}$, $ j=1,2, \ldots $,  be the coefficient matrices in the power series expansions of  the inverse matrix polynomial $ (I_m -\sum_{j=1}^p A_{j,p}z^j)^{-1}$, $ (I_m -\sum_{j=1}^p \widetilde{A}_{j,p}z^j)^{-1}$ and $ (I_m -\sum_{j=1}^p \widehat{A}_{j,p}z^j)^{-1}$, respectively,  $ |z|\leq 1$, where $ I_m $ is the $m\times m$ unit matrix.  
Set  $\Psi_{0,p}=\widetilde{\Psi}_{0,p}=\widehat{\Psi}_{0,p} =I_m$ and  let $ \Sigma_{\bm e}=E({\bm e}_{t,p}{\bm e}^\top_{t,p})$,  $ \widetilde{\Sigma}_{\bm e}=E(\widetilde{\bm e}_{t,p}
\widetilde{\bm e}^\top_{t,p})$ and $ \widehat{\Sigma}_{\bm e}=E(\widehat{\bm e}_{t,p}
\widehat{\bm e}^\top_{t,p})$. Since 
\[\Gamma_m(0)=\sum_{j=0}^\infty \Psi_{j,p} \Sigma_{\bm e}\Psi_{j,p}^\top, \  \ \Gamma^+_m(0)=\sum_{j=0}^\infty \widetilde{\Psi}_{j,p} \widetilde{\Sigma}_{\bm e}\widetilde{\Psi}_{j,p}^\top \ \ \mbox{and}  \ \ \Gamma^*_m(0)=\sum_{j=0}^\infty \widehat{\Psi}_{j,p} \widehat{\Sigma}_{\bm e}\widehat{\Psi}_{j,p}^\top,\]
the assertion of the lemma follows   using  the same  arguments as in the proof of  Lemma 6.5 of \cite{Paparoditis18}  and   the bounds 
\[\sum_{j=1}^\infty \| \widetilde{\Psi}_{j,p}-\Psi_{j,p}\|_F  = o_P\Big(m^{3/2}/\sqrt{p}\Big) \ \ \mbox{and} \ \  \sum_{j=1}^\infty \| \widehat{\Psi}_{j,p}-\widetilde \Psi_{j,p}\|_F= O_P\Big( \frac{p^5\sqrt{m}}{\sqrt{n}\lambda_m^2}\sqrt{\sum_{j=1}^m \alpha_j^{-2}}\Big),  \]
obtained in the aforementioned paper. \hfill $\Box$

\begin{lem} \label{le.gh} Let $ \widehat{g}, g  : {\mathcal H}^k \rightarrow {\mathcal H}$.
For  $ x=(x_1, x_2, \ldots, x_k) \in {\mathcal H}^k$ and $ h \in \mathbb{N}$, let
\[ \widehat{g}_{(h)}(x) = \widehat{g}\big(\widehat{g}_{(h-1)}(x), \widehat{g}_{(h-2)}(x), \ldots, (\widehat{g}_{(h-k)}(x) \big),\]
where $ \widehat{g}_{(1)}(x)=\widehat{g}(x)$ and $\widehat{g}_{(s)}(x)=x_{1-s}$ if $ s \leq 0$. Define analogously $ g_{(h)}(x)$. 
If
$ \|\widehat{g} - g_0\|_{\mathcal L} \stackrel{P}{\rightarrow} 0$ then,   
$ \|\widehat{g}_{(h)}(x) - g_{0,(h)}(x)\|_2 \stackrel{P}{\rightarrow} 0$  for any $ h \in \mathbb{N}$.
\end{lem}
{\bf Proof:} We use the notation  $ y_n=(\widehat{g}_{(L)}(x), \ldots, \widehat{g}_{(L-k+1)}(x)) $ and $ y= (g_{0,(L)}(x),  \ldots, g_{0,(L-k+1)}(x))$.  For $h=1$  we have 
\[ \|\widehat{g}(x) - g_{0}(x)\|_2 \leq \|\widehat{g}-g_0\|_{\mathcal L} \|x\|_2 = o_P(1)O(1) =o_P(1),\]
since $ \|\widehat{g}-g_0\|_{\mathcal L} \stackrel{P}{\rightarrow}0$. Suppose that $  \|\widehat{g}_{(L)}(x) - g_{0,(L)}(x)\|_2 \stackrel{P}{\rightarrow} 0$ for some $ L\in \mathbb{N}$. Then,
\begin{align*}
\| \widehat{g}_{(L+1)}(x) & - g_{0,(L+1)}(x)\|_{2}  \leq \|\widehat{g}(y_n) -  g_0(y_n) \|_2  \\
&  \ \ \ \  +  \|g_0(y_n) -  g_0(y) \|_2\\
& \leq \|\widehat{g} -g_0\|_{\mathcal L} \| y_n  \|_2 + \|g_0\|_{\mathcal L}\sum_{l=0}^{k-1}\| \widehat{g}_{(L-l)}(x) - g_{0,(L-l)}(x)\|_2 \\
& = o_P(1)O_P(1)  +  O(1)o_{P}(1) =o_{P}(1). 
\end{align*}
\hfill $\Box$

\noindent {\bf Proof of Proposition 3.1:} \ We first show that under the conditions of the proposition, and for every $ h\in {\mathbb N}$ the following  is true,
\begin{equation}\label{eq.prop1}  \sup_{\tau\in [0,1]}E \int_0^1\big(\widetilde{c}_{g,h}(\tau,s)- c_{g,h}(\tau,s)\big)^2ds \rightarrow 0,
\end{equation}
where $c_{g,h}$  and $ \widetilde{c}_{g,h}$ are the kernels associated with the operators $ g_{(h)} $and  $ \widetilde{g}_{(h)}$ respectively. These kernels   satisfy 
$c_{g,h}(\tau,s)= \int_0^1 c_{g,h-1}(\tau,u)c_{g,1}(u,s)du$ with $c_{g,1}(\tau,s)=c_g(\tau,s)$ and  $\widetilde{c}_{g,h}$  defined similarly.   Notice that for $ h=1$,
assertion~\eqref{eq.prop1} is true by assumption. Suppose that the assertion is true for some  $h \in {\mathbb N}$. Using  $ (a+b)^2 \leq 2a^2 + 2 b^2$ as well as Cauchy-Schwarz's inequality, we get
\begin{align*}
  \int_0^1\big(\widetilde{c}_{g,h+1}(\tau,s)- c_{g,h+1}(\tau,s)\big)^2ds  & \leq  2\int_0^1\big(\widetilde{c}_{g,h}(\tau,u)- c_{g,h}(\tau,u)\big)^2du  \int_0^1\int_0^1\widetilde{c}^2_{g}(u,s)duds \\
   &  + 2\int_0^1 c^2_{g,h}(\tau,u)du \int_0^1\int_0^1 \big(\widetilde{c}_{g}(u,s) -c_{g}(u,s)\big)^2duds.
\end{align*}  
From this we have 
\begin{align*}
 \sup_{\tau\in [0,1]}E \int_0^1\big(\widetilde{c}_{g,h+1}(\tau,s)- c_{g,h+1}(\tau,s)\big)^2ds  & \leq 2 \sup_{\tau\in [0,1]}E \int_0^1\big(\widetilde{c}_{g,h}(\tau,u)- c_{g,h}(\tau,u)\big)^2du\|\widehat{g}\|_{HS}^2\\
 &  + 2 \sup_{\tau\in [0,1]} \int_0^1c^2_{g,h}(\tau,u)du E\|\widetilde{g}-g\|^2_{HS} \\
  \leq  C \big\{\sup_{\tau\in [0,1]}E \int_0^1\big(\widetilde{c}_{g,h}(\tau,u) & - c_{g,h}(\tau,u)\big)^2du 
   +  E\|\widetilde{g}-g\|^2_{HS}\big\} \rightarrow 0.
 \end{align*}
Using~\eqref{eq.prop1} we  prove by induction that Assumption 3$^{'}$ is satisfied.  For $h=1$ we have
\begin{align*}
\sup_{\tau\in[0,1]}E\big|\widetilde{g}(x)(\tau) &- g(x)(\tau)\big|^2 \leq \|x\|_2^2 \sup_{\tau\in[0,1]}E\int_0^1\big(\widetilde{c}_g(\tau,s) -c_g(\tau,s)\big)^2ds  \rightarrow 0.
\end{align*}
Suppose that the assertion is true for some  $h\in {\mathbb N}$. Using Cauchy-Schwarz's inequality and the bound $ \|g(x)\|\leq \|g\|_{HS}\|x\|$, we  get for $ h+1$,
\begin{align*}
\sup_{\tau\in[0,1]}E\big|\widetilde{g}_{(h+1)}(x)(\tau)  & - g_{(h+1)}(x)(\tau)\big|^2  \leq  2\sup_{\tau\in[0,1]}E\big| \int_0^1\big(\widetilde{c}_{g,h}(\tau,s) -c_{g,h}(\tau,s)\big) \widetilde{g}(x)(s)ds \big|^2  \\
& \ \ \ + 2\sup_{\tau\in[0,1]}E\big| \int_0^1 c_{g,h}(\tau,s)\big( \widetilde{g}(x) -g(x)\big)(s)ds \big|^2 \\
& \leq  2\sup_{\tau\in[0,1]}E \int_0^1\big(\widetilde{c}_{g,h}(\tau,s) -c_{g,h}(\tau,s)\big)^2ds  \int_0^1\widetilde{g}^2(x)(s)ds \\
& \ \ \ + 2\sup_{\tau\in[0,1]}\int_0^1 c_{g,h}^2(\tau,s) E\int_0^1\big( \widetilde{g}(x) -g(x)\big)^2(s)ds \\
& \leq 2\|x\|_2^2 \sup_{\tau\in[0,1]}E \int_0^1\big(\widetilde{c}_{g,h}(\tau,s) -c_{g,h}(\tau,s)\big)^2ds \|\widetilde{g}\|^2_{HS}  \\
& \ \ \ \ + 2\|x\|_2^2\sup_{\tau\in[0,1]}\int_0^1 c_{g,h}^2(\tau,s) E\|\widetilde{g} -g\|^2_{HS} \rightarrow 0.
\end{align*}
 \hfill $\Box$

\noindent {\bf Proof of Theorem 3.1:} \  
Recall the notation  $ \X_{n,k}=(\X_{n-k+1}, \X_{n-k+2}, \ldots, \X_n)$. Observe that      
\begin{align*}
 \|\widehat{g}_{(h)}( {\X}_{n,k})-\widehat{g}^*_{(h)}( {\X}_{n,k}) \|_2  & \leq  \|\widehat{g}_{(h)}( {\X}_{n,k})-g_{0,(h)}( {\X}_{n,k}) \|_2 +
   \|\widehat{g}^*_{(h)}( {\X}_{n,k})-g_{0,(h)}( {\X}_{n,k}) \|_2 \\
 & =o_P(1),
\end{align*}
where the last equality follows by Assumption 3, Lemma~\ref{le.gh} and   the fact that $\|{\X}_{n,k}\|_2=O_P(1)$. Hence  
$
\mathcal E_{n+h}- \mathcal E_{n+h}^*  = \sum_{j=1}^\infty {\bf 1}^\top_{j}\bm{\xi}_{n+h}v_j 
-    \big(\sum_{j=1}^m {\bf 1}^\top_{j}\bm{\xi}^*_{n+h}\widehat{v}_j +U^*_{n+h,m}\big) +o_P(1)
$
 and by Slutsky's theorem, it suffices  to show that 
\begin{equation} \label{eq.d2-conv}
 d\big( \sum_{j=1}^\infty {\bf 1}_{j}^\top\bm{\xi}_{n+h}v_j   ,  \sum_{j=1}^m {\bf 1}_{j}^\top\bm{\xi}^*_{n+h}\widehat{v}_j + U^*_{n+h,m}\big)=o_P(1). 
 \end{equation}

 \noindent Assertion~\eqref{eq.d2-conv}  follows if  we  show that, 
\begin{enumerate}
\item[(i)] \ $ d\big(\sum_{j=1}^m \xi^+_{j,n+h}v_j  ,  \sum_{j=1}^\infty \xi_{j,n+h}v_j\big) \rightarrow 0$,
\item[(ii)]\  $ \|\sum_{j=1}^m {\bf 1}^\top_{j}\bm{\xi}^*_{n+h}\widehat{v}_j -\sum_{j=1}^m {\bf 1}^\top_{j}\bm{\xi}^+_{n+h}v_j \|_2 \stackrel{P}{\rightarrow}  0$, and,
\item[(iii)] \  $ U^*_{n+h, m} \stackrel{p}{\rightarrow}  0$.   
\end{enumerate}
To establish (i) consider  the sequence $ \{Y^+_{n}(h)\}$ in ${\mathcal H}$,  where 
$ Y^+_n(h)=\sum_{j=1}^\infty \widetilde{\xi}_{j,n+h}v_j$ and $ \widetilde{\xi}_{j,n+h} = \xi^+_{j,n+h}$ for  $ j=1,2, \ldots, m$ and  $ \widetilde{\xi}_{j,n+h}=0$ for $ j\geq m+1$.
For $ k \in \mathbb N$, let $ Y^+_{n,k}(h)= \sum_{j=1}^k \widetilde{\xi}_{j,n+h}v_j$. By  Theorem 3.2 of \cite{Billingsley99}, assertion (i) follows if we show that 
\begin{enumerate}
\item[(a)] \ $ Y^+_{n,k} (h)\stackrel{d}{\rightarrow} Y_k(h)=\sum_{j=1}^k \xi_{j,n+h}v_j$ for any $ k\in \mathbb N$, as $ n \rightarrow \infty$.
\item[(b)] \ $ Y_k(h) \stackrel{d}{\rightarrow} Y (h)=\sum_{j=1}^\infty \xi_{j,n+h} v_j$,  as $ k \rightarrow \infty$.
\item[(c)] \ For any $ \epsilon>0$, $ \lim_{
k\rightarrow \infty} \lim \sup_{n\in \mathbb N} P(\|Y^+_{n,m}(h)-Y_{n,k}^+(h) \| _2 >\epsilon) =0$.
\end{enumerate}
Consider (a). Assume that $n$ is large enough such that $m>k$. Since $k$ is fixed, $  (\widetilde{\xi}_{1,n+h}, \widetilde{\xi}_{2,n+h},$ $  \ldots, \widetilde{\xi}_{k,n+h})^\top $ $ =(\xi^+_{1,n+h}, \xi^+_{2,n+h}, \ldots, \xi^+_{k,n+h})^\top $ $ ={\bm \xi}^+_{n+h}(k)$ where the latter vector is obtained as
${\bm \xi}^+_{n+h}(k)= I_{k,m}{\bm \xi}^+_{n+h}$ with $ I_{k,m} $ the $k\times m$ matrix with  elements $(j,j)$, $j=1,2,\ldots, k$, equal to one and zero else.  Assume that $h=1$. Then, 
the  vector ${\bm \xi}^+_{n+1}$ is  
generated via the  regression type autoregression $ {\bm \xi}^+_{n+1} =  \sum_{j=1}^p\widetilde{A}_{j,p} \xi_{n+1-j} + e^+_{n+1}$ and  $ e^+_{n+1}$   i.i.d. innovations. Therefore,  and since $k$ is fixed, we have  by standard arguments (see  Lemma 3.1 of \cite{MK15}),  that $ {\bm \xi}^+_{n+1}(k) \stackrel{d}{\rightarrow} {\bm \xi}_{n+1}(k) =(\xi_{1,n+1}, \xi_{2,n+1}, \ldots, \xi_{k,n+1})^\top$. Suppose that the assertion is true for some $h \in {\mathbb N}$.  For $h+1$ it follows by similar arguments and using the recursion $ {\bm \xi}_{n+h+1}^+ = \sum_{j=1}^p \widetilde{A}_{j,p}{\bm \xi}^{+}_{n+h+1-j} + e^+_{n+h+1}$ that   
$ {\bm \xi}^+_{n+h+1}(k) \stackrel{d}{\rightarrow} {\bm \xi}_{n+h+1}(k) =(\xi_{1,n+h+1}, \xi_{2,n+h+1}, \ldots, \xi_{k,n+h+1})^\top$.
By the continuous mapping theorem we then conclude that $ Y^+_{n,k}(h) \stackrel{d}{\rightarrow} \sum_{j=1}^k\xi_{j,n+h}v_j =Y_{n,k}(h)$. 
Consider (b). Notice that, 
$ E\|Y^+_k (h)- Y(h)\|_2^2 = E\|\sum_{j=k+1}^\infty \xi_{j,n+h}v_j\|_2^2 = \sum_{j=k+1}^\infty \lambda_j \rightarrow 0$,
as $ k\rightarrow \infty$, which by Markov's inequality and Slusky's theorem implies that $ Y_k(h) \stackrel{d}{\rightarrow} \sum_{j=1}^\infty \xi_{j,n+h} v_j$ as $k\rightarrow \infty$.
Consider (c). 
We have 
\begin{align*}
E\|Y_{n,m}(h)-Y_{n,k}^+(h)\|^2_2 & = E\|\sum_{j=k+1}^m \xi^+_{j,n+h}v_j\|^2_2\\
& =\sum_{j=k+1}^m {\bf 1}_j^{\top} \Gamma^+_m(0) {\bf 1}_j  = \sum_{j=k+1}^m \lambda_j + \sum_{j=k+1}^m {\bf 1}_{j}^{\top}\big( \Gamma^+_m(0) -\Gamma_m(0)\big) {\bf 1}_j.
\end{align*}
Now since 
$  \big|\sum_{j=k+1}^m {\bf 1}_{j}^{\top}\big( \Gamma^+_m(0) -\Gamma_m(0)\big) {\bf 1}_j \big|= O_P\Big(\sqrt{m}\| \Gamma^+_m(0) -\Gamma_m(0)\|_F\big)$,
we get by Lemma~\ref{le.papa2018}, Assumption 2 and Markov's inequality that 
\[  \limsup\limits_{n\in \mathbb N} P(\|Y^+_{n,m}(h)-Y_{n,k}^+(h) \|_2>\epsilon) \leq \frac{1}{\epsilon^2} \big\{ \sum_{j=k+1}^\infty \lambda_j  + O_P\Big(\sqrt{m}\| \Gamma^+_m(0) -\Gamma_m(0)\|_F\big)\big\},\]
which converges to zero as $ k\rightarrow \infty$.

Consider assertion  (ii).  
 Since  $\|\widehat{v}_j\|_2=1$, we get the bound 
\begin{align*}
\|\sum_{j=1}^m {\bf 1}_{j}^\top\big(\bm{\xi}^+_{n+h}v_j -\bm{\xi}^*_{n+h}\widehat{v}_j\big)\|_2 &
 \leq \|\sum_{j=1}^m {\bf 1}_{j}^\top\big(\bm{\xi}^+_{n+h}-\bm{\xi}^*_{n+h}\big)\widehat{v}_j\|_2 \\
& \ \ \ \  + \|\sum_{j=1}^m {\bf 1}_{j}^\top\bm{\xi}^+_{n+h}\big(\widehat{v}_j-v_j\big)\|_2 \\
& \leq \sqrt{m}\|\bm{\xi}^*_{n+h}-\bm{\xi}^+_{n+h}\|_2 + \|{\bm \xi}^+_{n+h}\|_2\sum_{j=1}^m\|\widehat{v}_j-v_j\|_2\\
& =  \sqrt{m}\|\bm{\xi}^*_{n+h}-\bm{\xi}^+_{n+h}\|_2 + O_P\Big(\frac{m}{\sqrt{n}}\sqrt{\sum_{j=1}^m \alpha^{-2}_j}\Big),
\end{align*}
where  the last equality follows using $ \|{\bm \xi}^+_{n+h}\|_2^2 =O_P(m)$ and Lemma 3.2 of \cite{HK10}.
To evaluate the first term on the right hand side of the last displayed inequality, 
 we use the bound
\begin{align} \label{eq.decomE}
\|\bm{\xi}^*_{n+h}-\bm{\xi}^+_{n+h}\|_2 & \leq  \|\sum_{j=1}^p\big(\widehat{A}_{j,p}-\widetilde{A}_{j,p}\big){\bm{\xi}}^*_{n+h-j}\|_2 + \|\sum_{j=1}^p\widetilde{A}_{j,p}\big({\bm{\xi}}^+_{n+h-j} -\bm{\xi}_{n+h-j}^*\big)\|_2 
+ \|\bm{e}^*_{n+h}-\bm{e}^+_{n+h}\|_2\nonumber \\
& = \sum_{j=1}^3 T_{j,n},
\end{align}
with an obvious notation for $ T_{j,n}$, $ j=1,2,3$.  Observe first that 
 $ E\|\bm{e}^*_{n+h}-\bm{e}^+_{n+h}\|_2^2 \rightarrow 0$ 
in probability, as in 
the proof of Lemma 6.7 in \cite{Paparoditis18}, that is $ T_{3,n}\stackrel{P}{\rightarrow} 0$.
For the first term of~\eqref{eq.decomE} we have 
\begin{align*}
T_{1,n} &= \|\sum_{j=1}^p\big(\widehat{A}_{j,p}-\widetilde{A}_{j,p}\big){\bm{\xi}}^*_{n+h-j}\|_2  \leq \sum_{j=1}^p\|\widehat{A}_{j,p}-\widetilde{A}_{j,p}\|_F\|{\bm{\xi}}^*_{n+h-j}\|_2\\
& = O_P \big(\sqrt{mp}\|\widehat{A}_{p,m}-\widetilde{A}_{p,m}\|_F\big)
 = O_P\Big(\frac{m^{3/2} p^{5/2}}{\lambda_m^2} \sqrt{\frac{1}{n}\sum_{j=1}^m \alpha_{j}^{-2}}\Big) \rightarrow 0, 
\end{align*}
because 
$\|\widehat{A}_{p,m}-\widetilde{A}_{p,m}\|_F=O_P\big((p\sqrt{m}\lambda_m^{-1} +p^2)^2\sqrt{n^{-1}\sum_{j=1}^m \alpha_{j}^{-2}}\big),$ see \citet[][p. 5 of the Supplementary Material]{Paparoditis18}, and the last equality follows 
by Assumption 2.
To show that $ T_{2,n}\stackrel{P}{\rightarrow} 0$, we first show that, for any $h\in {\mathbb N}$,  
\begin{equation} \label{eq.T2.1}
\sum_{j=1}^p\|{\bm \xi}^+_{n+h-j}  -{\bm \xi}^*_{n+h-j}\|_2^2 \stackrel{P}{\rightarrow} 0,
\end{equation}
For $h=1$ we have 
\begin{align*}
\sum_{j=1}^p\|{\bm \xi}^+_{n+h-j}  -{\bm \xi}^*_{n+h-j}\|_2^2 & = \sum_{j=1}^p\|{\bm \xi}_{n+h-j}  -\widehat{{\bm \xi}}_{n+h-j} \|_2^2\\
& \leq \sum_{j=1}^p\|\X_{n+h-j}\|^2_2 \sum_{j=1}^m\|\widehat{v}_{j} - v_j  \|^2_2\\
& = O_P\big(p\sum_{j=1}^m \|\widehat{v}_j - v_j\|^2_2\big) = O_P\big(\frac{p}{n}\sum_{j=1}^m \alpha^{-2}_j \big) \rightarrow 0.
\end{align*}
Assume that  the assertion is true for some $h \in {\mathbb N}$. We then have for $h+1$ that
\begin{align*}
\sum_{j=1}^p\|{\bm \xi}^+_{n+h+1-j}  -{\bm \xi}^*_{n+h+1-j}\|_2^2& = \|{\bm \xi}^+_{n+h}  -{\bm \xi}^*_{n+h}\|_2^2+ o_P(1)\\
&  \leq 4 \|\sum_{j=1}^p (\widehat{A}_{j,p}- \widetilde{A}_{j,p}){\bm \xi}^*_{n+h-j} \|_2^2 + 4 \|\sum_{j=1}^p \widetilde{A}_{j,p} ({\bm \xi}_{n+h-j}^+ - {\bm \xi}_{n+h-j}^*)\|_2^2 \\
& \ \ \  + 2 \|{\bf e}^+_{h+h}-{\bf e}^\ast_{n+h}\|_2^2 + o_P(1), 
\end{align*}
which converges to zero in probability by the same arguments as those used in the proof that $ T_{1,n}$ and $T_{3,n}$ converge to zero in probability and because 
\[ \|\sum_{j=1}^p \widetilde{A}_{j,p} ({\bm \xi}_{n+h-j}^+ - {\bm \xi}_{n+h-j}^*)\|_2^2 \leq \sum_{j=1}^p\|\widetilde{A}_{j,p}\|_F^2 \sum_{j=1}^p\|{\bm \xi}_{n+h-j}^+ - {\bm \xi}_{n+h-j}^*\|_2^2 = O_P(1)\cdot o_P(1),\]
using the fact that $ \sum_{j=1}^p\|\widetilde{A}_{j,p}\|_F^2=O_P(1)$ uniformly in $ m$ and $p$; see the proof of Lemma 6.5 in \cite{Paparoditis18}.  Using~\eqref{eq.T2.1}, we get for the term $T_{2,n}$,
\begin{align*}
T^2_{2,n}= \|\sum_{j=1}^p \widetilde{A}_{j,p} ( {\bm{\xi}}^+_{n+h-j}- \bm{\xi}^*_{n+h-j})\|^2_2 & \leq \sum_{j=1}^p\|\widetilde{A}_{j,p}\|^2_F\sum_{j=1}^p \| {\bm{\xi}}^+_{n+h-j}- \bm{\xi}^*_{n+h-j}\|^2_2
 = O_P(1)\cdot o_P(1) \rightarrow 0.
\end{align*}

It remains to prove assertion (iii). This  
follows from Markov's inequality and  the fact that $E^\ast\|U^*_{n+h,m}\|^2_2 \rightarrow 0$ in probability. The last statement is true 
since 
\begin{align*}
E^\ast\|U^*_{n+h,m}\|^2_2& = \frac{1}{n}\sum_{t=1}^n\|\widehat{U}_{t,m} -\overline{U}_n\|_2^2\\
&  \leq  \frac{4}{n}\sum_{t=1}^n\|\widehat{U}_{t,m} -U_{t,m}\|_2^2  + \frac{4}{n}\sum_{t=1}^n\|U_{t,m}\|^2_2 + 2\|\overline{U}_n\|_2^2\\
&  = \frac{4}{n}\sum_{t=1}^n\|\sum_{j=1}^m {\bf 1}_{j}^\top\big(\bm{\xi}_{t}v_j -\widehat{\bm{\xi}}_{t}\widehat{v}_j\big)\|_2^2 + o_P(1) 
\end{align*}
where the $o_P(1)$ is due to  the weak law of large numbers,  the fact that $ E\|U_{t,m}\|^2 =\sum_{j=m+1}^\infty \lambda_j \rightarrow 0$ as $ m \rightarrow \infty$
and $ \overline{U}_{n,m} \rightarrow 0$, in probability.
For the first term  on the right hand side of the last displayed equality we have that this term is bounded by  
\begin{align*}
&  \frac{8}{n}\sum_{t=1}^n\big\|\sum_{j=1}^m {\bf 1}^\top_j(\widehat{\bm{\xi}}_t-\bm{\xi}_t)\widehat{v}_j\big\|_2^2+
\frac{8}{n}\sum_{t=1}^n\big\|\sum_{j=1}^m {\bf 1}^\top_j\bm{\xi}_t (\widehat{v}_j-v_j)\big\|_2^2 \\
&\leq \frac{8}{n}\sum_{t=1}^n\|\X_t\|^2_2\sum_{j=1}^m\|\widehat{v}_j-v_j\|_2^2   + 8\sqrt{\sum_{j,l=1}^m \big|\frac{1}{n}\sum_{t=1}^n\xi_{j,t}\xi_{l,t}\big|_2^2}
 \sum_{j=1}^m\|\widehat{v}_j-v_j\|_2^2 = O_P(n^{-1}\sum_{j=1}^m\alpha_j^{-2}) \rightarrow 0,
\end{align*} 
 where the last equality follows because $ n^{-1}\sum_{t=1}^n\|\X_t\|_2^2=O_P(1)$ and $n^{-1}\sum_{t=1}^n\xi_{j,t}\xi_{l,t} \stackrel{P}{\rightarrow} \lambda_j {\bf 1}_{j=l}$ as $ n \rightarrow \infty$.  \hfill $\Box$

\noindent  {\bf Proof of Theorem 3.2:} \  We only give the proof of assertion  (20).
 We have
 \begin{align} \label{eq.BasicBounSsigma}
 \big|\sigma_{n+h}^{*^2}(\tau) & -  \sigma_{n+h}^2(\tau)\big| \leq  \big| E^*(\X^*_{n+h}(\tau)^2 - E(\X_{n+h}(\tau)^2\big| \nonumber\\
 & \ \ + 2 \big| EX_{n+h}(\tau)\big(\widehat{g}_{(h)}(\X_{n,k})(\tau) - g_{0,(h)}(\X_{n,k})(\tau) \big) \big|
 \nonumber \\
 & \ \   + 2\big| E^*X_{n+h}^*(\tau)\big(\widehat{g}^*_{(h)}(\X_{n,k})(\tau) - g_{0,(h)}(\X_{n,k})(\tau) \big)\big| \nonumber \\
 &  + \big| E^*(\widehat{g}^*_{(h)}(\X_{n,k})(\tau) )^2- g^2_{0,(h)}(\X_{n,k})(\tau) \big| + \big| E(\widehat{g}_{(h)}(\X_{n,k})(\tau) )^2- g^2_{0,(h)}(\X_{n,k})(\tau) \big|
 \end{align}
 Notice first that  $ \sup_{\tau\in[0,1]} E^*(X_{n+h}^*(\tau))^2 =O_P(1)$ by equation~\eqref{eq.letzte} below  and the fact that  
 $   \sup_{\tau\in[0,1]} E(X_{n+h}(\tau))^2 =  \sup_{\tau \in [0,1]} c(\tau,\tau)  <\infty $, by the continuity of the kernel $c$. Using Caushy-Schwarz's inequality, and Assumption 3$^{'}$  we get, 
 \begin{align*}
 \sup_{\tau\in[0,1]}\big| EX_{n+h}(\tau) & \big(\widehat{g}_{(h)}(\X_{n,k})(\tau)  - g_{0,(h)}(\X_{n,k})(\tau) \big)\big|  \\
&  \leq  \sup_{\tau\in[0,1]} E(X_{n+h}(\tau))^2 \sqrt{\sup_{\tau\in [0,1]}
E\big( \widehat{g}_{(h)}(\X_{n,k})(\tau) - g_{0,(h)}(\X_{n,k})(\tau)\big)^2 } \rightarrow 0,
\end{align*} 
\begin{align*}
 \sup_{\tau\in[0,1]}\big| E^*X_{n+h}^*(\tau) & \big(\widehat{g}^*_{(h)}(\X_{n,k})(\tau)  - g_{0,(h)}(\X_{n,k})(\tau) \big)\big|  \\
&  \leq  \sup_{\tau\in[0,1]} E^*(X_{n+h}^*(\tau))^2 \sqrt{\sup_{\tau\in [0,1]}
E^*\big( \widehat{g}^*_{(h)}(\X_{n,k})(\tau) - g_{0,(h)}(\X_{n,k})(\tau)\big)^2 } \rightarrow 0.
\end{align*} 
Using $a^2-b^2=(a-b)(a+b)$, Cauchy-Schwarz's inequality  and Assumption 3$^{'}$  again, we get for  the last two terms on the right hand side of the bound~\eqref{eq.BasicBounSsigma}, that they also converges uniformly to zero, in probability.
 To establish assertion (20) it remains  to show that 
 \begin{equation} \label{eq.letzte}
  \sup_{\tau\in [0,1]} \big| E^*(\X^*_{n+h}(\tau))^2 - E(\X_{n+h}(\tau))^2\big|  \rightarrow 0,
  \end{equation}
   in probability.
Notice first that  due to the independence of $\sum_{j=1}^m {\bf 1}^\top {\bm \xi}^*_{n+h} \widehat{v}_j $ and $U^*_{n+h,m} $ we have 
$ E^*(\X^*_{n+h}(\tau))^2 = E^\ast( \sum_{j=1}^m {\bf 1}^\top {\bm \xi}^*_{n+h} \widehat{v}_j (\tau))^2 + E^*(U_{n+h,m}^*(\tau))^2$. Furthermore, using
$ c(\tau,\tau)=\sum_{j=1}^\infty \lambda_j v^2_j(\tau)$, where the convergence $|c(\tau,\tau)-\sum_{j=1}^k \lambda_j v^2_j(\tau)| \rightarrow 0$, as $ k \rightarrow \infty$,   is uniformly in $\tau\in[0,1]$, see Theorem 7.3.5 of \cite{HE15}, we get that 
 $ \sup_{\tau\in [0,1}E(U_{n+h,m}(\tau))^2 = \sup_{\tau\in [0,1]} \sum_{j=m+1}^\infty \lambda_j v^2_j(\tau) \rightarrow 0$.
 Therefore, and because  $ E(X_{n+h,m}(\tau) U_{n+h,m}(\tau))=0$ for all $\tau\in[0,1]$,   to establish~\eqref{eq.letzte}, it suffices by Cauchy-Schwarz inequality and  the inequality 
 $ \sup_{\tau\in[0,1]}\sqrt{f(\tau)}  \leq $ $   \sqrt{\sup_{\tau\in[0,1]}f(\tau)} $, where $f$ is a non negative  function on $[0,1]$, to show that, in probability, 
 \begin{enumerate}
 \item[(I)] \ $ \sup_{\tau\in [0,1]}  \big|E\big( \sum_{j=1}^m {\bf 1}_j^\top \big({\bm \xi}^*_{n+h} \widehat{v}_j (\tau)\big)^2 
 -E\big(\sum_{j=1}^m {\bf 1}_j^\top {\bm \xi}_{n+h} v_j (\tau)\big)^2\big| \rightarrow 0$, and 
 \item[(II)]  \ $ \sup_{\tau\in [0,1]} E^*\big(U_{n+h,m}^*(\tau)\big)^2 \rightarrow 0.$
 \end{enumerate}
Consider (I).  We have 
\begin{align*}
\big|E\big( \sum_{j=1}^m {\bf 1}_j^\top \big({\bm \xi}^*_{n+h} \widehat{v}_j (\tau)\big)^2 
 - & E\big(\sum_{j=1}^m {\bf 1}_j^\top {\bm \xi}_{n+h}   v_j (\tau)\big)^2\big| \leq   \big| \sum_{j_1,j_2=1}^m {\bf 1}^\top_{j_1}(\Gamma^*_m(0)-\Gamma_m(0)){\bf 1}_{j_2} v_{j_1}(\tau)v_{j_2}(\tau)\big|\\
 & \ \ \ \ +  \big| \sum_{j_1,j_2=1}^m {\bf 1}^\top_{j_1}\Gamma^*_m(0){\bf 1}_{j_2}\big( \widehat{v}_{j_1}(\tau)\widehat{v}_{j_2}(\tau)-  v_{j_1}(\tau)v_{j_2}(\tau)\big)\big|\\
 &  \leq \|\Gamma^*_m(0)-\Gamma_m(0)\|_F \big(\sum_{j=1}^m |v_j(\tau)| \big) \\
 &  \ \ \ \ +  \|\Gamma_m(0)\|_F\sum_{j=1}^m| \widehat{v}_j(\tau) -v_j(\tau)| \big( \sum_{j=1}^m |v_j(\tau)| + \sum_{j=1}^m |\widehat{v}_j(\tau)|\big).
\end{align*}
To evaluate the above terms notice that  $\|\Gamma_m(0)\|_F=O(1)$ and  that by  Lemma~\ref{le.papa2018}, $\|\Gamma^*_m(0)\|_F=O_P(1)$, where both  bounds are uniform in $m$. Furthermore,  using $ c(\tau,\tau) =\sum_{j=1}^\infty \lambda_j v^2_j(\tau)$ we get by the continuity of the kernel $c(\cdot,\cdot)$ on the compact support $ [0,1]\times[0,1]$, the bound
$ \sum_{j=1}^m \lambda_j v^2_j(\tau) \leq c(\tau,\tau)  \Rightarrow \sup_{\tau\in[0,1]}\sum_{j=1}^m v^2_j(\tau) \leq \lambda_m^{-1} C$, 
where $C := \sup_{\tau\in[0,1]}c(\tau,\tau) <\infty$. Moreover,  arguing as in \cite{Kokoszka2013}, 
we get that  $\{m^{-1}\sum_{j=1}^m (\sqrt{n}(\widehat{v}_j(\tau)-v_j(\tau)))^2, \tau\in [0,1]\}$ converges weakly which can be sharpened to  $ \sup_{\tau\in[0,1]}  m^{-1}\sum_{j=1}^m n(\widehat{v}_j(\tau)-v_j(\tau))^2 =O_P(1)$. 
Hence  
\[ \sum_{j=1}^m |\widehat{v}_j(\tau)| \leq \sqrt{m}\sqrt{\sum_{j=1}^m v^2_j(\tau)} + \sqrt{m} \sqrt{\sum_{j=1}^m (\widehat{v}_j(\tau) -v_j(\tau))^2} = 
O\Big( \sqrt{\frac{m}{\lambda_m}} + \frac{m}{\sqrt{n}}\Big),\]
and 
\[ \sum_{j=1}^m| \widehat{v}_j(\tau) -v_j(\tau)| \big( \sum_{j=1}^m |v_j(\tau)| + \sum_{j=1}^m |\widehat{v}_j(\tau)|\big) =  O_P\Big(\frac{m^{3/2}}{n \lambda^{1/2}_m } + \frac{m^2}{n^{3/2}} \Big).\]
Therefore,  
\[ \sup_{\tau\in[0,1]} \big|E\big( \sum_{j=1}^m {\bf 1}_j^\top \big({\bm \xi}^*_{n+h} \widehat{v}_j (\tau)\big)^2 
 -  E\big(\sum_{j=1}^m {\bf 1}_j^\top {\bm \xi}_{n+h}   v_j (\tau)\big)^2\big| = O_P\Big(\sqrt{\frac{m}{\lambda_m}} \|\Gamma^*_m(0)-\Gamma_m(0)\|_F\Big)  + o_P(1),\]
which vanishes because of  Lemma~\ref{le.papa2018} and Assumption 2$^{'}$.

Consider (II). Since  $ E^*(U^*_{n+h,m}(\tau))^2  \leq 2 n^{-1}\sum_{t=1}^n \big(\widehat{U}^2_{t,n} (\tau)\big)^2 + 2 \big(\overline{U}_n(\tau)\big)^2$, 
it suffices to show that $ \sup_{\tau\in[0,1]} n^{-1}\sum_{t=1}^n\big(\widehat{U}_{t,n}(\tau) \big)^2 \rightarrow 0$.  Toward this we use the bound 
\begin{align} \label{eq.decU}
\frac{1}{n}\sum_{t=1}^n\big(\widehat{U}_{t,n}(\tau) \big)^2 & \leq \frac{4}{n} \sum_{t=1}^n \|\widehat{\bm \xi}_t -{\bm \xi}_t \|^2_2 \sum_{j=1}^m v^2_j(\tau) + \frac{4}{n}\sum_{t=1}^n \|\widehat{\bm \xi}_{t}\|_2^2 \sum_{j=1}^m (\widehat{v}_j(\tau) - v_j(\tau))^2 \nonumber \\
& \ \  \  +  \frac{2}{n} \sum_{t=1}^n \big(\sum_{j=m+1}^\infty \xi_{j,t} v_j(\tau) \big)^2.
\end{align}
Since 
\[ \sup_{\tau\in[0,1]}\frac{1}{n} \sum_{t=1}^n \|\widehat{\bm \xi}_t -{\bm \xi}_t \|^2_2 \sum_{j=1}^m v^2_j(\tau) \leq C\frac{1}{\lambda_m}\frac{1}{n}\sum_{t=1}^n\|\widehat{\bm \xi}_t -{\bm \xi}_t \|^2_2 = \frac{1}{\lambda_m}O_P\Big(\frac{1}{n}\sum_{j=1}^m\frac{1}{\alpha^2_j}\Big) \]
and 
\[ \sup_{\tau\in[0,1]} \frac{1}{n}\sum_{t=1}^n \|\widehat{\bm \xi}_{t}\|_2^2 \sum_{j=1}^m (\widehat{v}_j(\tau) - v_j(\tau))^2 \leq O_P\Big(\frac{m}{n}\Big) \frac{1}{n}\sum_{t=1}^n\|\widehat{\bm \xi}_t\|_2^2 = O_P\Big(\frac{m}{n} \Big),\]
the first two terms on the right hand side of~\eqref{eq.decU} converge to zero. For the third term we get after evaluating the squared term and substituting 
 $ \xi_{j,t}=\langle \X_t, v_j\rangle$ the  bound, 
\begin{align}\label{eq.Bound-last}
\frac{1}{n}\sum_{t=1}^n(\X_t(\tau) & -\sum_{j=1}^m \xi_{j,t}v_j(\tau))^2  \leq \big|\frac{1}{n}\sum_{t=1}^n \X^2_t(\tau) -E\X_t^2(\tau)\big| \nonumber  \\
& \ \ + \big| \sum_{j_1=1}^m\sum_{j_2=1}^m \langle (\widehat{\mathcal C}_0-{\mathcal C}_0)(v_{j_1}), v_{j_2} \rangle v_{j_1}(\tau)v_{j_2}(\tau) \big| \nonumber \\
& \ \ + 2\big|\sum_{j=1}^m \int_0^1( \widehat{c}(\tau,s)-c(\tau,s))v_j(s)ds v_j(\tau) \big| + \big| E\X^2_t(\tau) - \sum_{j=1}^m \lambda_j v_j^2(\tau)\big|.
\end{align}
Now,  $ \sup_{\tau\in[0,1]}|n^{-1} \sum_{t=1}^n(\X^2_t(\tau) - E\X^2_t(\tau))| =o_P(1) $
and for the last term of~\eqref{eq.Bound-last} we have,  $ \sup_{\tau\in[0,1]}\big| E\X^2_t(\tau) - \sum_{j=1}^m \lambda_j v_j^2(\tau)\big| \rightarrow 0$, by Theorem 7.3.5 of \cite{HE15}. Also 
\begin{align*}
  \sup_{\tau\in[0,1]}\big| \sum_{j_1=1}^m\sum_{j_2=1}^m \langle (\widehat{\mathcal C}_0  -{\mathcal C}_0)(v_{j_1}), v_{j_2} \rangle v_{j_1}(\tau)v_{j_2}(\tau) \big|   
  & \leq \|\widehat{\mathcal C}_0-{\mathcal C}_0\|_{HS}\sup_{\tau\in[0,1]}\Big(\sum_{j=1}^m |v_j(\tau)|\Big)^2\\
  &  \leq  m \|\widehat{\mathcal C}_0-{\mathcal C}_0\|_{HS}\sup_{\tau\in[0,1]}\sum_{j=1}^m v^2_j(\tau) \\
  & \leq  C \frac{m}{\lambda_{m}} \|\widehat{\mathcal C}_0-{\mathcal C}_0\|_{HS}= O_{P}\Big(\frac{m}{\sqrt{n} \, \lambda_m} \Big),
  \end{align*}
 which converges to zero.  Finally, using
 \begin{align*}
 \big|\sum_{j=1}^m \int_0^1( \widehat{c}(\tau,s)-c(\tau,s))v_j(s)ds v_j(\tau) \big|^2 & \leq \int_0^1(\widehat{c}(\tau,s) -c(\tau,s))^2 ds \Big(\sum_{j=1}^m |v_j(\tau)| \Big)^2\\
 &  \leq \int_0^1(\widehat{c}(\tau,s) -c(\tau,s))^2 ds\, m \sum_{j=1}^m v^2_j(\tau),
 \end{align*}
 we get 
 \begin{align*}
 \sup_{\tau\in[0,1]}\big|\sum_{j=1}^m \int_0^1( \widehat{c}(\tau,s)-c(\tau,s))v_j(s)ds v_j(\tau) \big|  & \leq C\sqrt{\frac{m}{\lambda_m}} \sqrt{\sup_{\tau\in[0,1]}\int_0^1(\widehat{c}(\tau,s)-c(\tau,s))^2ds }\\
 & = O_P\Big(\sqrt{\frac{m}{\lambda_m\, n}}\Big) \rightarrow 0.
 \end{align*}
 \hfill $\Box$

 \noindent  {\bf Proof of Corollary 3.1:}  We  elaborate  only on  the case $h=1$.  In view of  (\ref{eq.unifvar}), the important step is   to show that $ {\mathcal E}_{n+1}^\ast$ converges weakly to the same limit as  $ {\mathcal E}_{n+1}$ in the space of continuous functions on  the interval $ [0,1]$, equipped with the supremum norm $\|\cdot\|_\infty$, i.e.,  $ (C([0,1], \|\cdot\|_\infty)$.  Toward this we assume that the random elements considered belong to this space.
 Since we do not consider any particular estimator $ \widehat{g}$ and $ \widehat{g}^\ast$, we assume that  $ g_0({\bf x})(\cdot)$ is uniformly continuous and that the corresponding sequences of estimators  fulfill the  property $ E \sup_{\tau\in[0,1]}\big|\widehat{g}({\bf x}) (\tau)-g_0({\bf x})(\tau)\big|  \rightarrow 0$, ${\bf x}  \in {\mathcal H}$, and the same  holds true for the estimator $ \widehat{g}^\ast$.  Notice that, this property has to be verified in a case by case 
 manner and depends on the functional $ g$ and the estimator  $\widehat{g}$ used. For instance, this assumption on $ \widehat{g}$ in the case of the estimator $\widehat{\varphi}_M(\tau,\sigma)$ discussed in Section 3.3, essentially can be derived from  an analogue assumption on the sequence of estimators  $\{\widehat{v}_j(\tau), \tau\in[0,1]\}$.  

The weak convergence of the finite dimensional distributions $ \big({\mathcal E}_{n+1}^\ast(\tau_j), j=1,2, \ldots, K\big)$ to the same limit as $  \big({\mathcal E}_{n+1}(\tau_j), j=1,2, \ldots, K\big)$,  for every positive integer $ K$, can  easily   be established by the Cram\'er-Wold device and a  slight modification of the arguments used in the proof of  assertion (\ref{eq.d2-conv})  established in the proof of  Theorem 3.1. That is by considering the  weak convergence in ${\mathbb R}$  for  linear combinations of the  (real valued) random variables $ \sum_{j=1}^m {\bf 1}_j^\top{\bm \xi}^\ast_{n+1}\widehat{v}_j(\tau_s)$ and $ \sum_{j=1}^\infty{\bf 1}_j^\top{\bm \xi}_{n+1}v_j(\tau_s)$,  $ s=1,2, \ldots, K$. To establish the equicontinuity condition,  that  for all $ \eta >0$ and all $ \epsilon >0$, there exists $\delta >0$, such that 
    \[ \limsup_{n\rightarrow\infty} P^\ast(\sup_{|\tau_1-\tau_2| <  \delta} |{\mathcal E}_{n+1}^\ast(\tau_1)-{\mathcal E}_{n+1}^\ast (\tau_2)| >\epsilon) <\eta,\]
 use the decomposition 
   \begin{align} \label{eq.E-dec}
 {\mathcal E}_{n+1}^\ast(\tau_1)  - {\mathcal E}^\ast_{ n+1}(\tau_2) &  = \sum_{j=1}^m\xi_{j,t}^\ast \big(\widehat{v}_j(\tau_1) -   \widehat{v}_j(\tau_2) \big) \nonumber  \\
 & \ \ \ \ +  \big(\widehat{g}^\ast({\bf X}_{n,k})-\widehat{g}^\ast({\bf X}_{n,k})(\tau_2)\big) +  \big(U^\ast_{n+1}(\tau_1) - U^\ast_{n+1}(\tau_2)\big).
 \end{align}
   By Markov's inequality, the  first term can be bounded by 
 \begin{align*}
 E^\ast\big(\sup_{|\tau_1-\tau_2| <  \delta} \big|\sum_{j=1}^m\xi_{j,t}^\ast \big(\widehat{v}_j(\tau_1) -   \widehat{v}_j(\tau_2) \big)\big| )  & \leq 2 \sum_{j=1}^m \sup_{\tau \in [0,1]}|\widehat{v}_j(\tau) -v_j(\tau)|\big(E^\ast(\xi_{j,t}^\ast)^2\big)^{1/2}  \\
 & \ \ \ +  \sum_{j=1}^m \sup_{|\tau_1-\tau_2| < \delta} |v_{j}(\tau_1)-v_j(\tau_2)|\big(E^\ast(\xi_{j,t}^\ast)^2\big)^{1/2} , 
 \end{align*}
 which can be made arbitrarily small  in probability,  using the arguments of the proof of Theorem 3.2,  the  consistency   properties of the $\widehat{v}_j(\cdot)$'s and the fact that $E^\ast(\xi_{j,t}^\ast)^2\rightarrow 0$, in probability, as $ j \rightarrow \infty$.   
 The second term  in (\ref{eq.E-dec}) is bounded by,  (recall ${\bf X}_{n,k} $ is fixed), 
\begin{align*}
E^\ast\big(\sup_{|\tau_1-\tau_2| <  \delta} \big|\widehat{g}^\ast({\bf X}_{n,k}) (\tau_1)-\widehat{g}^\ast({\bf X}_{n,k})(\tau_2)\big|\big) & \leq 2 E^\ast  \sup_{\tau\in[0,1]}\big|\widehat{g}^\ast({\bf X}_{n,k}) (\tau)-g({\bf X}_{n,k})(\tau)\big| \\
& \ \ \ \ +  \sup_{|\tau_1-\tau_2| < \delta}|g({\bf X}_{n,k})(\tau_1)-g({\bf X}_{n,k})(\tau_2)|.
\end{align*}
In   bounding   the  third  term in  decomposition (\ref{eq.E-dec}), it turns out that the corresponding bound    is dominated  by the term  
$ n^{-1}\sum_{t=1}^n \sup_{|\tau_1-\tau_2| < \delta} \big|U_{t,m}(\tau_1)-U_{t,m}(\tau_2)\big| \leq   2n^{-1}\sum_{t=1}^n \|U_{t,m}\|_\infty $, which converges to zero in probability, by  the 
fact that $ \|U_{t,m}\|_\infty  \rightarrow 0$ in probability, as  $m\rightarrow \infty$. 
%
%

\bibliographystyle{agsm}
\bibliography{sieve}

\end{document}